\documentclass[preprint,sort&compress,12pt]{elsarticle}
\usepackage{amsmath}
\usepackage{mathrsfs}
\usepackage{stmaryrd}
\usepackage{mathrsfs}
\usepackage{bm}
\usepackage{amsmath}
\usepackage{amsfonts}
\usepackage{amssymb}
\usepackage{amsthm}

\journal{\quad}
\newcommand{\mf}{\mathbf}
\newcommand{\mm}{\mathrm}

\begin{document}
\begin{frontmatter}
\title{On
Instability and Stability of Three-Dimensional Gravity \\
Driven Viscous Flows in a Bounded Domain
}

\author[FJ]{Fei Jiang\corref{cor1}}
\ead{jiangfei0591@163.com}
\cortext[cor1]{Corresponding author: Tel +86-18305950592.}
\author[sJ]{Song Jiang}
\ead{jiang@iapcm.ac.cn}
\address[FJ]{College of Mathematics and Computer Science, Fuzhou University, Fuzhou, 350108, China.}
\address[sJ]{Institute of Applied Physics and Computational Mathematics, 
 Beijing, 100088, China.}
\begin{abstract}
We investigate the instability and stability of some steady-states 
of a three-dimensional nonhomogeneous incompressible viscous flow driven by gravity
in a bounded domain $\Omega$ of class $C^2$. When the steady density is
heavier with increasing height (i.e., the Rayleigh-Taylor steady-state),
we show that the steady-state is linear unstable (i.e., the linear solution grows in
time in $H^2$) by constructing a (standard) energy functional and exploiting the modified variational
method. Then, by introducing a new energy functional and using a careful bootstrap argument,
we further show that the steady-state is nonlinear unstable in the sense of Hadamard.
When the steady density is lighter with increasing height, we show, with the help of a restricted
condition imposed on steady density, that the steady-state is linearly globally stable and nonlinearly
locally stable in the sense of Hadamard.
\end{abstract}

\end{frontmatter}


\newtheorem{thm}{Theorem}[section]
\newtheorem{lem}{Lemma}[section]
\newtheorem{pro}{Proposition}[section]
\newtheorem{cor}{Corollary}[section]
\newproof{pf}{Proof}
\newdefinition{rem}{Remark}[section]
\newtheorem{definition}{Definition}[section]

\section{Introduction}
\label{Intro} \numberwithin{equation}{section}

The motion of a three-dimensional (3D) nonhomogeneous incompressible
viscous fluid in the presence of a uniform gravitational field in a
bounded domain $\Omega\subset {\mathbb R}^3$ of $C^2$-class is governed by the following
Navier-Stokes equations:
\begin{equation}\label{0101}\left\{\begin{array}{l}
 \rho_t+{\bf v}\cdot\nabla \rho=0,\\[1mm]
\rho\mathbf{v}_t+\rho {\bf v}\cdot\nabla {\bf v}+\nabla p=\mu\Delta\mathbf{v}-g\rho\mf{e}_3,\\[1mm]
\mathrm{div}\mathbf{v}=0,\end{array}\right.\end{equation}
where the unknowns $\rho:=\rho(t,\mf{x})$, $\mathbf{v}:=\mathbf{v}(t,\mf{x})$ and $p:=p(t,\mf{x})$ denote the density,
velocity and pressure of the fluid, respectively;
$\mu>0$ stands for the coefficient of shear viscosity, $g>0$ for the gravitational
constant, $\mf{e}_3=(0,0,1)$ for the vertical unit vector, and $-g\mf{e}_3$ for the gravitational force.
In the system (\ref{0101}) the equation (\ref{0101})$_1$ is the continuity equation, while
\eqref{0101}$_2$ describes the balance law of momentum.

The stability/instability of viscous incompressible flows governed by the Navier-Stokes
equations is a classical subject with a very extensive literature over more than 100 years.
In this paper we study the instability and stability of the following steady-state to the system \eqref{0101}:
\begin{eqnarray*}
&&\mathbf{v}(t,\mathbf{x})\equiv \mathbf{0},\\
&&\nabla \bar{p}=-\bar{\rho}g \mf{e}_3\quad \mbox{ in }\Omega.
\end{eqnarray*}
It is easy to show that the steady density $\bar{\rho}$ only depends on $x_3$, the third component of $\mf{x}$,
provided $\bar{\rho}\in C^1(\Omega)$. Hence we denote
$\bar{\rho}':=\partial_{x_3}\bar{\rho}$ in this paper for simplicity.
Moreover, we can compute out the corresponding steady pressure $\bar{p}$ determined by $\bar{\rho}$.
Now, denote the perturbation by
$$ \varrho=\rho -\bar{\rho},\quad \mathbf{u}=\mathbf{v}-\mathbf{0},\quad q=p-\bar{p},$$
then, $(\varrho ,\mathbf{u},q)$ satisfies the perturbed equations:
\begin{equation}\label{0105}\left\{\begin{array}{l}
\varrho_t+{\bf u}\cdot\nabla (\varrho+\bar{\rho})=0, \\[1mm]
(\varrho+\bar{\rho}){\bf u}_t+(\varrho+\bar{\rho}){\bf u}\cdot\nabla
{\bf u}+\nabla q=\mu\Delta{\mathbf{{u}}}-g\varrho\mf{e}_3,\\[1mm]
 \mathrm{div}\mathbf{u}=0.\end{array}\right.  \end{equation}
To complete the statement of the perturbed problem, we specify the
initial and boundary conditions:
\begin{equation}\label{0106}
(\varrho,{\bf u})|_{t=0}=(\varrho_0,{\bf u}_0)\quad\mbox{in } \Omega
\end{equation}
and
\begin{equation}\label{0107}
{\bf u}|_{\partial\Omega}={\bf 0}\quad \mbox{ for any }t>0.
\end{equation}
Moreover, the initial data should satisfy the compatibility
condition $\mathrm{div}\mathbf{u}_0=0$.

If we linearize the equations (\ref{0105}) around the steady state
$(\bar{\rho},\mathbf{0})$, then the resulting linearized equations read as
\begin{equation}\label{0108}
\left\{\begin{array}{ll}
 \varrho_t+ \bar{\rho}'u_3=0, \\[1mm]
  \bar{\rho}\mathbf{u}_t +\nabla q=\mu\Delta{\mathbf{{u}}}-g\varrho\mf{e}_3,\\[1mm]
  \mathrm{div}\mathbf{u}=0,
\end{array}\right.\end{equation}
where $u_3$ denotes the third component of $\mf{u}$.

First we consider the case that the steady-state is unstable, this situation corresponds to
the well-known Rayleigh-Taylor (RT) instability. So, we consider the following RT density profile
$\bar{\rho}:=\bar{\rho}(\mf{x})$ which exists and satisfies
\begin{eqnarray}\label{0102}
\bar{\rho}\in C^{2}(\bar{\Omega}),\ \inf_{\mf{x}\in
 \Omega}\{\bar{\rho}(\mf{x})\}>0,\end{eqnarray}
 and
\begin{eqnarray}\label{0103}\bar{\rho}'(\mf{x}_0)>0\quad\mbox{ for some }\mf{x}_0\in \Omega.
\end{eqnarray}
The condition (\ref{0103}) means that there is a region in which the RT density profile has
larger density with increasing $x_3$ (height), thus leading to the linear and nonlinear
RT instability as shown in Theorem \ref{thm:0101} and \ref{thm:0102} below.

The RT instability is well known as gravity-driven instability in fluids
when a heavy fluid is on top of a light one. Instability of the linearized problem (i.e. linear instability) for an
incompressible fluid was first introduced by Rayleigh in 1883
\cite{RLAP}. In 2003, Hwang and Guo \cite{HHJGY} proved the nonlinear RT
 instability of $\|(\varrho,\mf{u})\|_{L^2(\Omega)}$ in the sense of Hadamard
 for a 2D nonhomogeneous incompressible inviscid fluid (i.e. $\mu=0$ in (\ref{0105})--(\ref{0106}))
 with boundary condition $\mathbf{u}\cdot \mathbf{n}|_{\partial\Omega}=0$, where
$\Omega=\{(x_1,x_2)\in \mathbb{R}^2~|~-l<x_2<m\}$ and $\mathbf{n}$ denotes the
 outer normal vector to $\partial\Omega$. Recently, Jiang et.al. \cite{NJTSC2} showed the nonlinear RT
 instability of $\| {u}_3\|_{L^2(\mathbb{R}^3)}$ for the Cauchy problem of (\ref{0105})
 in the sense of Lipschitz structure, and further gave the nonlinear RT
 instability of $\|u_3\|_{L^2(\Omega)}$ in \cite{JFJSWWWN} in the sense of Hadamard in a
unbounded horizontal period domain $\Omega$.
In addition, similar results of linear or nonlinear RT instability were established for two
layer incompressible or compressible fluids with a free interface (so-called stratified fluids),
where the RT steady-state solution is a denser fluid lying
above a lighter one separated by a free interface and the domain is also a flat domain (such as $\mathbb{R}^3$,
an infinite slab domain and a horizontal period domain), please refer to
\cite{PJSGOI5,GYTI1,wang2011viscous,GYTI2} for incompressible and compressible stratified fluids,
and \cite{KMSMSP,HRWP,JFJSWWWOA,DRJFJS} for stratified MHD fluids.

To our best knowledge, however, it is still open mathematically whether there exists an unstable
solution to the linearized problem of (\ref{0105})--(\ref{0107}) and the nonlinear problem
(\ref{0105})--(\ref{0107}) of a 3D viscous fluid in a general bounded domain.
One of the aims in the present paper is to show rigorously
the instability for the initial-boundary problem (\ref{0105})--(\ref{0107}) in a bounded domain
in the sense of Hadamard. There are a number of articles studying the instability of flows in a flat domain,
see, for example, the instability for the periodic BGK equilibria \cite{GYSWIC},
for the horizontal period Rayleigh-B\'enard convection \cite{GYHYQCQ}, and
for the space periodic quasi-geostrophic equation \cite{FSNPVVNC},
for an ideal space periodic fluid \cite{FSSWVMNA,EG200,VMFSN2003} and for
the space periodic and whole space forced MHD equations \cite{BIIS413,GDOS2}.
Compared with the flat domain case,
essential difficulties in techniques arise for general domains, because one can not apply the method of the
Fourier transform (or discrete mode-$e^{i\xi\cdot\mf{x}}$) to analyze properties of spectrums of the
associated linearized problem. Here we shall modify and adapt the modified variational method \cite{GYTI2} for ODEs to our
partial differential equations problem
to circumvent such difficulties in the study the linear RT instability.

Before stating our main results, we explain the notations
used throughout this paper. For simplicity, we
drop the domain $\Omega$ in Sobolve spaces and the corresponding norms as well as in
integrands over $\Omega$, for example,
\begin{equation*}  \begin{aligned}&
L^p:=L^p(\Omega),\quad
{H}^1_0:=W^{1,2}_0(\Omega),\;\;
{H}^k:=W^{k,2}(\Omega),\;\; {H}^{\sigma}_0:=\{\tilde{\mathbf{v}}\in {H}^1_0(\Omega)~|~
\mm{div}\tilde{\mathbf{v}}=0\},\;\; \int:=\int_\Omega .  \end{aligned}\end{equation*}
In addition, a product space $(X)^n$ of vector functions
are still denoted by $X$, for examples, the vector function $\mf{u}\in (H^2)^3$ is denoted
by $\mf{u}\in H^2$ with norm $\|\mf{u}\|_{H^2}:=(\sum_{k=1}^3\|u_k\|_{H^2}^2)^{1/2}$.

 The first main result is on the linear instability and reads as follows.
\begin{thm}[Linear instability]\label{thm:0101}
 Assume that the steady density profile $\bar{\rho}$ satisfies (\ref{0102}), (\ref{0103}).
Then the steady state $(\bar{\rho},\mathbf{0})$ of the linearized system (\ref{0106})--(\ref{0108})
 is linearly unstable. That is, there exists an unstable solution
$$(\mf{\varrho},\mf{u},{q}):=e^{\Lambda t}(\tilde{v}_3/\Lambda,\tilde{\mf{v}},\tilde{p})$$
 to \eqref{0106}--\eqref{0108}, where $(\tilde{\mf{v}},\tilde{p})\in H^2\times H^1$
 solves the following boundary problem
 \begin{equation}\label{0113}
\left\{
                              \begin{array}{ll}
\Lambda^2\bar{\rho}\tilde{\mathbf{v}}+\Lambda\nabla \tilde{{q}}
=\Lambda\mu\Delta\tilde{\mathbf{{v}}}+g
\bar{\rho}'\tilde{v}_3 \mf{e}_3,\\[1mm]
\mathrm{div}\, \tilde{\mathbf{v}}=0,\quad
\tilde{\mathbf{v}}|_{\partial\Omega}=0
\end{array}
                            \right.
\end{equation} with the constant growth rate $\Lambda$ defined by
\begin{equation}\begin{aligned}
\label{0111nn} \Lambda^2=\sup_{\tilde{\mf{w}}\in
{H}^\sigma_0(\Omega)}\frac{g \int \bar{\rho}'\tilde{{w}}_3^2\mm{d}\mf{x}-\Lambda\mu
\int|\nabla
\tilde{\mathbf{w}}|^2\mm{d}\mf{x}}{\int\bar{\rho}|\tilde{\mathbf{w}}|^2\mathrm{d}\mf{x}}
.
\end{aligned}\end{equation}
 Moreover, $\tilde{\mf{v}}$ satisfies
 $\tilde{{v}}_3\equiv\!\!\!\!\!\!/\ 0$, $
\tilde{{v}}_1^2+\tilde{{v}}_2^2\equiv\!\!\!\!\!\!/\ 0$. In particular,
\begin{equation} \|(\mf{\varrho},\mf{u},q)(t)\|_{L^2}\to \infty,\qquad\mbox{as }t\to\infty . \label{js1}
\end{equation}
\end{thm}
\begin{rem}
For any $\epsilon >0$, $(\epsilon\mf{\varrho},\epsilon\mf{u},\epsilon q)$ is still a solution of (\ref{0106})--(\ref{0108})
with the initial data that can be made arbitrarily small. This solution is obviously unstable thanks to
(\ref{js1}).
\end{rem}

We sketch the proof framework for Theorem \ref{thm:0101}.
To construct a solution to the
linearized equations (\ref{0108}) that has growing $H^k$-norm for
any $k$, we shall make a growing normal mode ansatz of solutions, i.e.,
\begin{equation*}
{\varrho}(\mathbf{x})=\tilde{\rho} (\mathbf{x})e^{\Lambda t},\;\;
\mathbf{u}(\mathbf{x})=\tilde{\mathbf{v}}(\mathbf{x})e^{\Lambda
t},\;\; {q}(\mathbf{x})=\tilde{p}(\mathbf{x})e^{\Lambda
t}\quad\mbox{for some }\Lambda>0.
\end{equation*}
Substituting this ansatz into (\ref{0108}), one obtains the time-independent system
\begin{equation}\label{0109}
\left\{
                              \begin{array}{ll}
\Lambda\tilde{\rho}+\bar{\rho}'\tilde{{v}}_3=0,\\[1mm]
\Lambda \bar{\rho}\tilde{\mf{v}}+ \nabla \tilde{{p}}
=\mu\Delta\tilde{\mathbf{{v}}}-g
\tilde{\rho}\mf{e}_3,\\[1mm]
\mathrm{div}\, \tilde{\mathbf{v}}=0,\quad \tilde{\mathbf{v}}|_{\partial\Omega}=\mathbf{0}.
\end{array}
                            \right.
\end{equation}
Then eliminating $\tilde{\rho}$ by using the first equation, we arrive at the boundary problem \eqref{0113}.
For the domain $\Omega$, which is not a flat domain, we can not adapt the methods as in \cite{GYTI2,GYTI1,NJTSC2,HHJGY},
where the boundary problem \eqref{0113} can be changed to ordinary
differential equations (ODEs) by employing the horizontal Fourier transform.
In \cite{GYTI2}, Guo and Tice used a modified variational method to construct a solution
to a class of ODEs arising from the linearized partial differential equations (PDEs).
In this paper, we extend the modified variational method for ODEs in \cite{GYTI2} to the PDEs \eqref{0113} to construct
a solution of \eqref{0113}. In view of the basic idea of the modified variational method,
we first modify the boundary problem \eqref{0109} as follows.
 \begin{equation}\label{nnn0113}
\left\{
                              \begin{array}{ll}
\Lambda^2\bar{\rho}\tilde{\mathbf{v}}+\Lambda\nabla \tilde{{p}}
=s\mu\Delta\tilde{\mathbf{{v}}}+g
\bar{\rho}'\tilde{v}_3 \mf{e}_3,\\[1mm]
\mathrm{div}\, \tilde{\mathbf{v}}=0,\quad
\tilde{\mathbf{v}}|_{\partial\Omega}=\mf{0},
\end{array}
                            \right.
\end{equation}
which has a standard energy $E(\tilde{\mathbf{v}},s)$ (see \eqref{0204}) and an associated
admissible set $\mathcal{A}$ (see \eqref{0205}). Denoting $\alpha(s):=\sup_{\tilde{\mathbf{v}}\in\mathcal{A}}E(\tilde{\mathbf{v}},s)$,
then by virtue of the variational method, there exists a pair of functions $(\tilde{\mathbf{v}},\tilde{p})$ which solves \eqref{nnn0113}
with $\Lambda=\sqrt{\alpha(s)}$ for each $s>0$. By investigating the properties of ${\alpha(s)}$,
we further find that there is a finite interval $(0,\mathfrak{S})$ such that ${\alpha(s)}\in C(0,\mathfrak{S})$ and
 $\alpha (s)>0$ on $(0,\mathfrak{S})$. Moreover, $\lim_{s\rightarrow 0}\alpha(s)>0$
 and $\lim_{s\rightarrow \mathfrak{S}}\alpha(s)=0$. It is not difficult to see that by a fixed point argument, we
 immediately get functions $(\tilde{\mathbf{v}},\tilde{p})$ satisfying \eqref{nnn0113} with $\Lambda=\sqrt{\alpha(\Lambda)}>0$.
 Thus we conclude Theorem \ref{thm:0101}.

Based on Theorem \ref{thm:0101}, we can establish the following nonlinear instability result.\begin{thm}[Nonlinear instability]\label{thm:0102}
 Assume that the RT density profile $\bar{\rho}$ satisfies (\ref{0102}) and
 \begin{eqnarray}\label{nnnn0103}\inf_{\mf{x}\in \Omega}\{\bar{\rho}'(\mf{x})\} >0.
\end{eqnarray}
 Then, the steady state $(\bar{\rho},\mathbf{0})$ of the system (\ref{0105})--(\ref{0107})
 is unstable in the Hadamard sense, that is, there are positive constants $\Lambda$, $m_0$, $\varepsilon$ and $\delta_0$,
 and functions $(\bar{\varrho}_0,\bar{\mf{u}}_0)\in H^2\times H^2$,
such that for any $\delta\in (0,\delta_0)$ and the initial data
 $(\varrho_0,\mathbf{u}_0):=(\delta\bar{\varrho}_0,\delta\bar{\mf{u}}_0)$
there is a unique strong solution $({\varrho},\mathbf{u})\in C^0([0,T^{\max}],L^2\times H^2)$ of (\ref{0105})--(\ref{0107})
with a associated pressure $p\in H^1$, such that
\begin{equation*}\label{0115}\|\varrho(T^\delta)\|_{L^2},\ \|(u_1,{u}_2)(T^\delta)\|_{L^2},\ \|{u}_3(T^\delta)\|_{L^2}\geq {\varepsilon}\;
\end{equation*}
for some escape time $T^\delta:=\frac{1}{\Lambda}\mm{ln}\frac{2\varepsilon}{m_0\delta}\in
(0,T^{\max})$, where $T^{\max}$ denotes the maximal time of existence of the solution
$(\varrho,\mathbf{u})$.
\end{thm}

{Now we sketch the main idea in the proof of Theorem \ref{thm:0101}. Similarly to \cite{GYTI2,HHJGY}, we can use
energy estimates to show that $\Lambda$ in $e^{\Lambda t}$ constructed in Theorem \ref{thm:0101} is the sharp growth rate in the
norm $\sqrt{\|\varrho\|_{L^2 }^2+\|\mathbf{u}\|_{H^2 }^2}$. This means that the initial velocity $\mathbf{u}_0\in H^2$ should be
divergence-free (i.e., the compatibility condition $\mm{div}{\mathbf{u}_0}=0$).
So, if we use Duhamel's principle in the standard bootstrap argument, the nonlinear term
``$( \varrho+\bar{\rho})\mathbf{u}\cdot\nabla\mathbf{u} +\varrho\mathbf{u}_t$'' in \eqref{0105} should
belong to $H^2$ and satisfy the compatibility condition of
divergence-free. Unfortunately, the nonlinear term does not possess the property of divergence-free in general. We remark here that
for weak solutions, this will not yields any difficulties because the compatibility conditions of high order are not needed as
studied in other articles (see \cite{BIIS413,GDOS2} for examples). To circumvent this obstacle,
we shall use some specific energy estimates to replace Duhamel's principle. However, in the derivation of the energy estimates,
we find that we have to control $\mf{u}_{tt}$, hence the nonlinear solution must be a classical solution, this requires that the initial velocity
of the nonlinear problem has to satisfy the higher compatibility conditions (such as ${\bf u}_t|_{t=0}=[-{\bf u}\cdot\nabla
{\bf u}+(\mu\Delta{\mathbf{{u}}}-g\varrho\mf{e}_3-\nabla q)/(\varrho+\bar{\rho})]|_{t=0}\in H_0^1(\Omega)$). However, the
initial velocity of the linear solution constructed in Theorem \ref{thm:0101}
could not satisfy the higher compatibility condition of the initial velocity for the nonlinear system \eqref{0105}--\eqref{0107}.
To overcome this difficulty, we introduce a new energy functional $E_{\mm{N}}(\tilde{\rho},\tilde{\mf{v}})$ (see \eqref{jiangfeijiangads})
with a new admissible set $\mathcal{A}_\mm{N}$ (see \eqref{n0205}) for the time-independet system
\eqref{0109} under the condition \eqref{nnnn0103}.
It is interesting that the growth rate re-defined by
\begin{equation}\label{newenergy1n}\Lambda:=\sup_{(\tilde{\rho},\tilde{\mathbf{v}})\in
\mathcal{A}_{\mm{N}}}E_\mm{N}(\tilde{\rho},\tilde{\mathbf{v}})
\end{equation}
is equal to the growth rate constructed by \eqref{0111nn} (see Proposition \ref{equal}).
 This fact is very important, since we shall further show that $\Lambda$ of $e^{\Lambda t}$ is the sharp growth rate in the
lower norm ``$\sqrt{\|\varrho\|_{L^2 }^2+\|\mathbf{u}\|_{L^2 }^2}$'' by \eqref{newenergy1n}.
This shows that the initial data of the linear solutions can be used as
the initial data of the strong solutions to the nonlinear system \eqref{0105}--\eqref{0107}.
Consequently, we can circumvent the problem of the higher compatibility condition. }

An interesting question arises whether the nonlinear system is stable when the steady density is
lighter with increasing height (i.e., $\bar{\rho}'(\mf{x})\leq 0$ for any $\mf{x}\in \Omega$).
Unfortunately, we still can not give a complete answer to this question. We have, however,
the following stability results by employing the standard energy method,
if $\bar{\rho}$ satisfies an additional condition.
\begin{thm}[Linear and nonlinear stability]\label{thm:0103}
Assume that $\bar{\rho}(\mf{x})$ satisfy \eqref{0102} and
  \begin{eqnarray}\label{supedenstiy}
 \sup_{\mf{x}\in \Omega}\{\bar{\rho}'(\mf{x})\}<0.\end{eqnarray}
\begin{enumerate}
  \item[(1)] Then the linearized system \eqref{0106}--\eqref{0108} is globally stable.
    That is, for any $(\varrho_0,\mathbf{u}_0)\in H^1\times H^2$, there exits a unique global strong solution
 $(\varrho,\mathbf{u})\in C^{0}((0,\infty),L^2\times H^2)$ to \eqref{0106}--\eqref{0108},
 such that for any $t>0$,
  \begin{eqnarray}\label{nnn0117}
&&\|(\varrho, \mathbf{u})(t)\|^2_{L^2}+\int_0^t\|\nabla
 \mathbf{u}\|^2_{L^2}\mm{d}s\leq C\|(\varrho_0, \mathbf{u}_0)\|^2_{L^2},\\
&&\label{nnn0118}\|\nabla\mf{u}(t)\|^2_{H^1}+\int_0^t\|\nabla\mathbf{u}_t  \|^2_{L^2}\mm{d}s
\leq C\|(\varrho_0,\mathbf{u}_0,\Delta\mf{u}_0)\|^2_{L^2},  \end{eqnarray}
where the constant $C$ only depends on $g$, $\mu$ and $\bar{\rho}$.
Moreover, $\mf{u}$ is asymptotically stable, i.e.,
 \begin{equation}\label{nnn0119}\|\mathbf{u}(t)\|_{H^1}\to 0\quad \mbox{ as }t\to\infty.\end{equation}
  \item[(2)] If we further assume that $\bar{\rho}'$ is constant, then the nonlinear system
\eqref{0105}--\eqref{0107} is locally stable.  That is, for any $K>0$, there is a $\delta_0:=\delta_0(K)>0$,
such that for any $(\varrho_0,\mathbf{u}_0)\in H^1\times H^2$ satisfying
 $$\|(\varrho_0,\mathbf{u}_0,\Delta\mf{u}_0)\|_{L^2}\leq \delta_0\;\;\mbox{ and }\; 0<\inf_{\mf{x}\in \Omega}\{(\varrho_0
 +\bar{\rho})(\mf{x})\}\leq \sup_{\mf{x}\in \Omega}\{(\varrho_0+\bar{\rho})(\mf{x})\}\leq K, $$
 there exists a unique global strong solution
 $(\varrho,\mathbf{u})\in C^{0}((0,\infty),L^2\times H^2)$ to \eqref{0105}--\eqref{0107} satisfying \eqref{nnn0119},
\eqref{nnn0117} with a different constant $C$ depending on $g$, $\mu$ and $\bar{\rho}$ only, and
  \begin{eqnarray}\label{131103fu}
\|\nabla\mf{u}(t)\|^2_{H^1}+\int_0^t\|\nabla\mathbf{u}_t\|^2_{L^2}\mm{d}s\leq C(K)\|(\varrho_0,
\mathbf{u}_0,\Delta\mf{u}_0)\|^2_{L^2},\qquad \forall\; t>0, \end{eqnarray}
 where the constant $\delta_0(K)$ depends on $g$, $\mu$ and $\bar{\rho}$, and is nonincreasing on $K$, and the constant $C(K)$
 depends on $g$, $\mu$ and $\bar{\rho}$, and is nondecreasing on $K$.
\end{enumerate}
\end{thm}
\begin{rem}
We point out that the smallness condition on the initial data in (2) of Theorem \ref{thm:0103} is needed
in the derivation of the higher order stable estimate \eqref{131103fu}.
\end{rem}

The rest of this paper is organized as follows: In Sections \ref{sec:02}, \ref{sec:05} and \ref{sec:06},
we give the proof of Theorem \ref{thm:0101}, \ref{thm:0102} and \ref{thm:0103}, respectively. In
Section \ref{sec:03} and \ref{sec:04}, we analyze the growth rates
and deduce the nonlinear energy estimates, which will be used in the proof of Theorem \ref{thm:0103}.

\section{Proof of the linear instability}\label{sec:02}
In this section, we adapt the modified variational method to show Theorem \ref{thm:0101}.
   First, multiplying \eqref{0113} by $\tilde{\mathbf{v}}$ and integrating the resulting identity, we get
\begin{equation*} \begin{aligned}
                            \Lambda^2\int
                            \bar{\rho}\tilde{\mathbf{v}}^2\mm{d}\mf{x}=g\int
\bar{\rho}'\tilde{v}_3^2\mm{d}\mf{x}
- s\mu \int|\nabla\tilde{\mf{v}}|^2\mm{d}\mf{x}.
\end{aligned}\end{equation*}
Hence the standard energy functional of (\ref{nnn0113}) can be given by
\begin{equation}\label{0204}E(\tilde{\mf{v}})=g\int \bar{\rho}'\tilde{v}_3^2\mm{d}\mf{x}-s\mu
\int|\nabla\tilde{\mf{v}}|^2\mm{d}\mf{x}\end{equation} with an associated admissible set
\begin{equation}\label{0205}
 \mathcal{A}:=\left\{\tilde{\mathbf{v}}\in
H^\sigma_0~\bigg|~J(\tilde{\mf{v}}):=\int
                            \bar{\rho}\tilde{\mathbf{v}}^2\mm{d}\mf{x}=1\right\}.
\end{equation}
Thus we can find a $\Lambda$ by maximizing
\begin{equation}\label{0206}
\Lambda^2:=\sup_{\tilde{\mathbf{v}}\in
\mathcal{A}}E(\tilde{\mathbf{v}}).\end{equation}
Obviously, $\sup_{\tilde{\mathbf{v}}\in\mathcal{A}}E(\tilde{\mathbf{v}})<\infty$ for any $s\geq 0$. In order to emphasize the dependence of
$E(\tilde{\mathbf{v}})$ upon $s>0$, we will sometimes write
\begin{equation*}E(\tilde{\mf{v}},s):=E(\tilde{\mf{v}})\quad\mbox{ and }\quad
\alpha(s):=\sup_{\tilde{\mathbf{v}}\in
\mathcal{A}}E(\tilde{\mathbf{v}},s). \end{equation*}

Next we show that a maximizer of (\ref{0206}) exists, and that the corresponding
Euler-Lagrange equations are equivalent to
(\ref{nnn0113}).
\begin{pro}\label{pro:0201}
 Assume that the steady density profile $\bar{\rho}$ satisfy
 \begin{eqnarray}\label{lsd0102}
\bar{\rho}\in C^{1}(\bar{\Omega}),\ \inf_{\mf{x}\in
 \Omega}\{\bar{\rho}(\mf{x})\}>0,\end{eqnarray}
then for any but fixed $s>0$, the following assertions are valid.
 \begin{enumerate}
    \item[(1)] $E({\tilde{\mathbf{v}}})$
achieves its supremum on $\mathcal{A}$.
   \item[(2)] Let $\tilde{\mathbf{v}}_0$ be
a maximizer and $\Lambda:=\sqrt{\sup_{\tilde{\mathbf{v}}\in\mathcal{A}}E(\tilde{\mathbf{v}})}>0$,
then there exists a corresponding pressure field $\tilde{p}_0$ associated to $\tilde{\mathbf{v}}_0$,
 such that the triple ($\tilde{\mathbf{v}}_0$, $\tilde{p}_0$, $\Lambda$) satisfies
the boundary problem (\ref{nnn0113}). Moreover,
$(\tilde{\mathbf{v}}_0,\tilde{p}_0) \in H^2\times H^1$.
 \end{enumerate}
\end{pro}
\begin{pf}
(1) 
Let $\tilde{\mathbf{v}}_n\in \mathcal{A}$ be a maximizing sequence, then
$E(\tilde{\mathbf{v}}_n)$ is lower bound. This together with
(\ref{0205}) implies that
$\tilde{\mathbf{v}}_n$ is bounded in $H^1$. So, there exists a
$\tilde{\mathbf{v}}_0\in H^1\cap\mathcal{A}$ and a subsequence (still denoted by $\mf{v}_n$ for simplicity), such that
$\tilde{\mathbf{v}}_n\rightarrow \tilde{\mathbf{v}}_0$ weakly in
$H^1$ and strongly in $L^2$. Moreover, by the lower semi-continuity, one has
\begin{equation*}
\begin{aligned}
\sup_{\tilde{\mathbf{v}}\in \mathcal{A}}E(\tilde{\mathbf{v}})= \limsup_{n\rightarrow
\infty}E(\tilde{\mathbf{v}}_n)
= &\lim_{n\rightarrow \infty}\int \bar{\rho}'\tilde{v}_{n3}^2\mm{d}\mf{x}-
s\mu\liminf_{n\rightarrow \infty} \int|\nabla\tilde{\mf{v}}_n|^2\mm{d}\mf{x}\\
\leq & E(\tilde{\mathbf{v}}_0)\leq \sup_{\tilde{\mathbf{v}}\in\mathcal{A}}E(\tilde{\mathbf{v}}),
\end{aligned}\end{equation*}
which shows that $E(\tilde{\mathbf{v}})$ achieves its supremum on $\mathcal{A}$.

(2) To show the second assertion, we notice that since $E(\tilde{\mathbf{v}})$ and $J(\tilde{\mathbf{v}})$
 are homogeneous of degree $2$, (\ref{0206}) is equivalent to
  \begin{equation}\label{0227}
  \Lambda^2=\sup_{\tilde{\mathbf{v}}\in {H}^{\sigma}_0}\frac{E(\tilde{\mathbf{v}})}{J(\tilde{\mathbf{v}})}.
\end{equation}
For any $\tau\in \mathbb{R}$ and $\mf{w}\in {H}^\sigma_0$, we take
$\tilde{\mf{w}}(\tau):=\tilde{\mathbf{v}}_0+\tau\mf{w}$. Then (\ref{0227}) implies
  \begin{equation*}E(\tilde{\mf{w}}(\tau))-\Lambda^2J(\tilde{\mf{w}}(\tau))\leq 0.
\end{equation*}
If we set
$I(\tau)=E(\tilde{\mf{w}}(\tau))-\Lambda^2J(\tilde{\mf{w}}(\tau))$,
then we see that $I(\tau)\in C^1(\mathbb{R})$, $I(\tau)\leq 0$ for all $\tau\in \mathbb{R}$ and
$I(0)=0$. This implies $I'(0)=0$. Hence, a direct computation leads to
  \begin{equation*}\begin{aligned}
\Lambda^2\int
                            \bar{\rho}\tilde{\mathbf{v}}_0\cdot\mf{w}\mm{d}\mf{x}
                            +s\mu\int \nabla\tilde{\mf{v}}_0:\nabla\mf{w}
\mm{d}\mf{x}=g\int
\bar{\rho}'\tilde{{v}}_{03}\mf{e}_3\cdot\tilde{\mf{w} }\mm{d}\mf{x}\end{aligned}  \end{equation*}
which implies that $\tilde{\mathbf{v}}$ is a weak solution to the boundary problem \eqref{nnn0113}.
Here $\tilde{{v}}_{03}$ denotes the third component of $\tilde{\mathbf{v}}_0$.
Thus, it follows from the classical regularity theory on the Stokes equations
(see \cite[Theorem IV.6.1]{galdi2011introduction}) and the condition \eqref{lsd0102} that
there are constants $c_1$ dependent of the domain $\Omega$), and $c_2$ dependent of $g$, $\Lambda$, $\Omega$
and $\bar{\rho}$, such that
  \begin{equation*}    
  \|\tilde{\mathbf{v}}\|_{H^{2}} +\Lambda\|\nabla \tilde{p}_0\|_{L^{2}}\leq c_1\|
(g\bar{\rho}'\tilde{{v}}_{03}\mf{e}_3-\Lambda^2\bar{\rho}\tilde{\mathbf{v}}_0)/(s\mu )\|_{L^2}\leq c_2,
\end{equation*}
  where $\tilde{p}_0 \in H^{1}$ is the corresponding pressure field associated to $\tilde{\mathbf{v}}_0$.
 This completes the proof.
\hfill $\Box$
\end{pf}

Next, we want to show that there is a fixed point such that $\Lambda=s>0$. To this end,
we first give some properties of $\alpha(s)$ as a function of $s> 0$.
\begin{pro}\label{pro:0202}
 Assume that the steady density profile $\bar{\rho}$ satisfies (\ref{lsd0102}). Then the
function $\alpha(s)$ defined on $(0,\infty)$ enjoys the following properties:
\begin{enumerate}[\quad \ (1)]
 \item $\alpha(s)\in C_{\mathrm{loc}}^{0,1}(0,\infty)$ is nonincreasing.
   \item
 There are constants $c_3$, $c_4>0$ which depend on $g$, $\bar{\rho}$ and $\mu$, such that
  \begin{equation}\label{0210}\alpha(s)\geq c_3-sc_4
  .\end{equation}
  \end{enumerate}
\end{pro}
\begin{pf}
(1) Let $\{\tilde{\mathbf{v}}^n_{s_i}\}\subset\mathcal{A}$ be a maximizing sequence
of $\sup_{\tilde{\mathbf{v}}\in\mathcal{A}}E(\tilde{\mathbf{v}},s_i)=\alpha(s_i)$ for $i=1$ and $2$. Then
\begin{equation*}
\alpha(s_1)\geq
\limsup_{n\rightarrow\infty}E(\tilde{\mathbf{v}}_{s_2}^n,s_1) \geq
\liminf_{n\rightarrow\infty}E(\tilde{\mathbf{v}}_{s_2}^n,s_2)=\alpha(s_2)\;
\mbox{ for any }0<s_1<s_2<\infty.
\end{equation*}
Hence $\alpha(s)$ is nonincreasing on $(0,\infty)$. Next we use this fact to show the continuity of $\alpha(s)$.

Let $I:=[a,b]\subset (0,\infty)$ be a bounded interval. In view of
\eqref{lsd0102} and the monotonicity of $\alpha(s)$, we know that
\begin{equation} \label{02321}
|\alpha(s)|\leq \max\left\{|\alpha(b)|,{g}  \left\|\bar{\rho}'/{\bar{\rho}}
\right\|_{L^\infty}\right\}<\infty\quad\mbox{ for any }s\in I.
\end{equation}
On the other hand, for any $s\in I$, there exists a maximizing sequence
$\{\tilde{\mathbf{v}}^n_{s}\}\subset\mathcal{A}$ of $\sup_{\tilde{\mathbf{v}}\in
\mathcal{A}}E(\tilde{\mathbf{v}},s)$, such that
\begin{equation}\label{0232}\begin{aligned}|\alpha(s)-E(\tilde{\mathbf{v}}_{s}^n,s)|<1
\end{aligned}.\end{equation}
Making use of (\ref{0204}), (\ref{02321}) and (\ref{0232}), we infer that
\begin{equation*}\begin{aligned}\label{0234}0\leq
\mu \int|\nabla \tilde{\mathbf{v}}|^2\mm{d}\mf{x}
=&\frac{g }{s}\int \bar{\rho}'|\tilde{ {v}}_{s3}^n|^2\mathrm{d}\mf{x} -\frac{E(\tilde{\mathbf{v}}_s^n,s)}{s}
\\
\leq &
\frac{1+\max\{|\alpha(b)|,{g}\left\|{\bar{\rho}'}/{\bar{\rho}}
\right\|_{L^\infty}\}}{a}+\frac{g}{a}\left\|\frac{\bar{\rho}'}{\bar{\rho}}
\right\|_{L^\infty }:=K.
\end{aligned}\end{equation*}
Thus, for $s_i\in I$ ($i=1,2$), we further find that
\begin{equation}\begin{aligned}\label{0235}\alpha(s_1)= \limsup_{n\rightarrow
\infty}E(\tilde{\mathbf{v}}_{s_1}^n,s_1)\leq & \limsup_{n\rightarrow
\infty}E(\tilde{\mathbf{v}}_{s_1}^n,s_2)+\mu
|s_1-s_2|\limsup_{n\rightarrow
\infty}\int_{\Omega}|\nabla \tilde{\mathbf{v}}_{s_1}^n|^2\mm{d}\mf{x}\\
\leq & \alpha(s_2)+K|s_1-s_2|.
\end{aligned}\end{equation}
Reversing the role of the indices $1$ and $2$ in the derivation of the inequality
(\ref{0235}), we obtain the same boundedness with the indices switched. Therefore, we deduce that
\begin{equation*}\begin{aligned}|\alpha(s_1)-\alpha(s_2)|\leq K|s_1-s_2|,
\end{aligned}\end{equation*}
which yields $\alpha(s)\in C_{\mathrm{loc}}^{0,1}(0,\infty)$.

(2) We turn to prove \eqref{0210}. We first construct a function ${\mf{v}}\in \mathcal{A}$ such that
\begin{equation}\label{0214}
g \int\bar{\rho}'{ {v}}_{3}^2 \mm{d}\mf{x}>0.
\end{equation}
Noting that since $\bar{\rho}'(\mf{x}_0)>0$ for some point $\mf{x}_0\in \Omega$ and the domain $\Omega$
is an open set, there is a ball $B_{\mf{x}_0}^{\delta}:=\{\mf{x}\in\Omega~|~|\mf{x}-\mf{x}_0|<\delta\}$,
such that
\begin{equation*}
\bar{\rho}'(\mf{x})>0\quad\mbox{ for any }\mf{x}\in B_{\mf{x}_0}^{\delta}.\end{equation*}
Now, choose a smooth function $f(r)\in C^1(\mathbb{R})$, such that
\begin{equation*}f(r)=-f(-r),\  |f(r)|>0\hbox{ if }0<|r|<{\delta}/{4},\
\mbox{ and }f(r)=0\mbox{ if }|r|\geq{\delta}/{4},
 \end{equation*}
and define then
\begin{equation*}\bar{\mf{v}}(\mf{x}):=f(x_1)\left(0,-
f(x_3)\int_{-{\delta/4}}^{x_2} f(r
)\mm{d}r,f(x_2)\int_{-{\delta/4}}^{x_3}f(r)\mm{d}r\right).\end{equation*}
It is easy to check that the non-zero function $\bar{\mf{v}}(\mf{x})\in H_0^\sigma(B_{\mf{0}}^{\delta})$, thus
${ {\mathbf{v}}}:=\mathbf{u}(\mf{x}-\mf{x}_0)\in
H_0^\sigma(\Omega)
$ and satisfies \eqref{0214}

With \eqref{0214} to hand, one has
\begin{equation*}\begin{aligned}
\alpha(s)=& \sup_{\tilde{\mathbf{v}}\in
\mathcal{A}}E(\tilde{\mathbf{v}},s)=\sup_{\tilde{\mathbf{v}}\in
{H}^\sigma_0}\frac{E(\tilde{\mathbf{v}},s)}{J(\tilde{\mathbf{v}})}
\\
&\geq \frac{E({  {\mathbf{v}}},s)}{J({ {\mathbf{v}}})}=
\frac{g  \int{\bar{\rho}'
{v}}_{3}^2\mm{d}\mf{x}}{\int
                            \bar{\rho}{\mathbf{v}}^2\mm{d}\mf{x}}
-s\frac{ \mu \int|\nabla
{\mathbf{v}}|^2\mm{d}\mf{x}}{\int
                            \bar{\rho}{\mathbf{v}}^2\mm{d}\mf{x}}:= c_3-sc_4
\end{aligned}\end{equation*} for two positive constants $c_3:
 =c_3(g,\bar{\rho})$ and $c_4:=c_4(g,\mu,\bar{\rho})$.
This completes the proof of Proposition \ref{pro:0202}.
 \hfill $\Box$
\end{pf}

Next we show that there exists a pair of functions $(\tilde{{v}},\tilde{p})$ satisfying
\eqref{0113} with grow rate $\Lambda$.  Let
\begin{equation}\label{0240}\mathfrak{S} :=\sup\{s~|~\alpha(\tau)>0\mbox{ for any }\tau\in (0,s)\}.
\end{equation}
By virtue of Proposition \ref{pro:0202}, $\mathfrak{S}>0$, and moreover, $\alpha(s)>0$ for any $s<\mathfrak{S}$.
Since $\alpha(s)=\sup_{\tilde{\mathbf{v}}\in\mathcal{A}}E(\tilde{\mathbf{v}},s)<\infty$,
using the monotonicity of $\alpha(s)$, we see that
 \begin{equation}\label{zero}
 \lim_{s\rightarrow 0}\alpha(s)\mbox{ exists and the limit is a positve constant.}
 \end{equation}

 On the other hand, by virtue of Poincar\'{e}'s inequality, there is a constant $c_5$ dependent of
$g$, $\bar{\rho}$ and $\Omega$, such that
 $$g\int \bar{\rho}'\tilde{v}_3^2\mm{d}\mf{x}\leq \left|g\int \bar{\rho}'\tilde{v}_3^2\mm{d}\mf{x}\right|\leq c_5
\int|\nabla\tilde{\mf{v}}|^2\mm{d}\mf{x}\quad\mbox{ for any }\tilde{\mathbf{v}}\in\mathcal{A}.$$
 Thus, if $s>c_5/\mu$, then
 $$g\int \bar{\rho}'\tilde{v}_3^2\mm{d}\mf{x}-s\mu
\int|\nabla\tilde{\mf{v}}|^2\mm{d}\mf{x}<0\quad\mbox{ for any }\tilde{\mathbf{v}}\in\mathcal{A},$$
which implies that
 $$\alpha(s)\leq 0\quad \mbox{ for any } s>c_5/\mu. $$
Hence $\mathfrak{S}<\infty$, moreover,
\begin{equation}\label{zerolin}
\lim_{s\rightarrow \mathfrak{S}}\alpha(s)=0.
\end{equation}
Now, using a fixed-point argument, exploiting \eqref{zero}, \eqref{zerolin}, and the continuity of $\alpha(s)$ on
$(0,\mathfrak{S})$, we find that there exists a unique $\Lambda\in(0,\mathfrak{S})$, such that
\begin{equation} \label{growth}
 \Lambda=\sqrt{\alpha(\Lambda)}=\sqrt{\sup_{\tilde{\mathbf{w}}\in
\mathcal{A}}E(\tilde{\mathbf{w}}, \Lambda)}>0. \end{equation}

In view of Proposition \ref{pro:0201}, there is a solution $(\tilde{\mf{v}},\tilde{p})\in H^2\times H^1$
to the problem (\ref{nnn0113}) with $\Lambda$ constructed in \eqref{growth}. Moreover,
$\Lambda^2=E(\tilde{\mf{v}},\Lambda)$,
$\tilde{{v}}_3\not\equiv 0$ by \eqref{growth} and \eqref{0204}. By the embedding theorem, $ \tilde{\mf{v}}\in C^{0,\varsigma}(\Omega)$ for some constant $\varsigma\in (0,1)$; moreover,
from the fact that $\mm{div}\tilde{\mathbf{v}}=0$
and $\tilde{\mathbf{v}}|_{\Omega}=\mf{0}$ we get $\tilde{{v}}_1^2+\tilde{{v}}_2^2\not\equiv 0$ immediately.
Thus we conclude the following proposition, which implies Theorem \ref{thm:0101}.
\begin{pro}\label{pro:nnn0203}
 Assume that the steady density profile $\bar{\rho}$ satisfies (\ref{lsd0102}). Then there exists a pair of functions
$(\tilde{\mf{v}},\tilde{p})\in H^{2}\times H^{1}$ satisfying the following boundary problem:
\begin{equation*}
\left\{
                              \begin{array}{ll}
\Lambda^2\bar{\rho}\tilde{\mathbf{v}}+\Lambda\nabla \tilde{{q}}
=\Lambda\mu\Delta\tilde{\mathbf{{v}}}+g
\bar{\rho}'\tilde{v}_3 \mf{e}_3,\\[1mm]
\mathrm{div}\, \tilde{\mathbf{v}}=0,\quad
\tilde{\mathbf{v}}|_{\partial\Omega}=0
\end{array}
                            \right.
\end{equation*}
with a growth rate $\Lambda>0$ defined by $\Lambda^2=\sup_{\tilde{\mathbf{w}}\in\mathcal{A}}E(\tilde{\mathbf{w}},\Lambda)$.
 Moreover, $\tilde{\mf{v}}$ satisfies $\tilde{{v}}_3\not\equiv 0$ and $\tilde{{v}}_1^2+\tilde{{v}}_2^2\not\equiv 0$.
 In particular, if $\bar{\rho}\in C^2(\bar{\Omega})$, then $\tilde{\rho}:=\bar{\rho}'\tilde{{v}}_3/\Lambda\in H^1$
 and $(\tilde{\rho},\tilde{\mf{v}},\tilde{p})$ satisfies
 \begin{equation*}    \begin{array}{ll}
\Lambda\tilde{\rho}+\bar{\rho}'\tilde{{v}}_3=0,\\[1mm]
\Lambda \bar{\rho}\tilde{\mf{v}}+ \nabla \tilde{{q}}=\mu\Delta\tilde{\mathbf{{v}}}
-g \tilde{\rho}\mf{e}_3.
\end{array}  \end{equation*}
\end{pro}

\section{Additional results on the sharp growth rate}\label{sec:03}
In this section, we shall show that $\Lambda$ constructed  by
(\ref{growth}) is the sharp growth rate of arbitrary solutions to the
linearized problem \eqref{0106}--(\ref{0108}). Since the density varies, the
spectrum of the linear operator is difficult to analyze, and it is hard to
obtain the largest growth rate of the solution operator in some Sobolev
space in the usual way. Instead, motivated by \cite{GYTI2}, we use
energy estimates to show that $e^{\Lambda t}$ is the sharp growth rate in
the norm $\sqrt{\|\varrho\|_{L^2 }^2+\|\mathbf{u}\|_{H^2 }^2}$.
Moreover, if $\bar{\rho}$ satisfies \eqref{nnnn0103}, we can further show that
$e^{\Lambda t}$ is the sharp growth rate in the lower norm $\|(\varrho,\mathbf{u})\|_{L^2 }$.

To obtain the sharp growth rate in the lower norm, we shall look for a new definition of $\Lambda$.
First we change the boundary problem \eqref{0109} into the following form:
 \begin{equation}\label{n0201}
\left\{
                              \begin{array}{ll}
\Lambda\tilde{\rho}/\bar{\rho}'=-\tilde{{v}}_3,\\[1mm]
(\Lambda \bar{\rho}\tilde{\mf{v}}+ \nabla \tilde{{q}}-\mu\Delta\tilde{\mathbf{{v}}})/g
=-
\tilde{\rho}\mf{e}_3,\\[1mm]
\mathrm{div}\, \tilde{\mathbf{v}}=0,\quad \tilde{\mathbf{v}}|_{\partial\Omega}=\mathbf{0}.
\end{array}
                            \right.
\end{equation}
 The new form \eqref{n0201} has a fine variational
structure under the condition \eqref{nnnn0103}. In fact, multiplying \eqref{n0201}$_1$ and  \eqref{n0201}$_2$  by $\tilde{\rho}$ and
$\tilde{\mathbf{v}}$, respectively, integrating the resulting identities,  and adding them, we get
\begin{equation*} \begin{aligned}
                            \Lambda\int
\left( \frac{\tilde{\rho}^2}{\bar{\rho}'}+\frac{\bar{\rho}|\tilde{\mathbf{v}}|^2}{g}\right)\mm{d}\mf{x}=
-\int\left(\frac{\mu |\nabla\tilde{\mf{v}}|^2}{g}+2\tilde{\rho}\tilde{v}_3\right)\mm{d}\mf{x}.
\end{aligned}\end{equation*}
So, we define an energy functional of
(\ref{n0201}) by
\begin{equation}\label{jiangfeijiangads}E_{\mm{N}}(\tilde{\rho},\tilde{\mf{v}})=
-\int\left(\frac{\mu |\nabla\tilde{\mf{v}}|^2}{g}+2\tilde{\rho}\tilde{v}_3\right)\mm{d}\mf{x}
\end{equation}
with an associated admissible set
\begin{equation}\label{n0205}
 \mathcal{A}_{\mm{N}}:=\left\{(\tilde{\rho}, \tilde{\mathbf{v}})\in
L^2\times H^{\sigma}_0~\bigg|~J_{\mm{N}}(\tilde{\rho},\tilde{\mf{v}}):=\int
\left( \frac{\tilde{\rho}^2}{\bar{\rho}'}+\frac{\bar{\rho}|\tilde{\mathbf{v}}|^2}{g}\right)\mm{d}\mf{x}=1\right\}.
\end{equation}
Thus we can find a $\Lambda$ by maximizing
\begin{equation}\label{n0206}
\Lambda:=\sup_{(\tilde{\rho},\tilde{\mathbf{v}})\in
\mathcal{A}_{\mm{N}}}E_\mm{N}(\tilde{\rho},\tilde{\mathbf{v}}).\end{equation}

By virtue of the condition of steady density $\bar{\rho}$,  obviously
$\sup_{(\tilde{\rho},\tilde{\mathbf{v}})\in\mathcal{A}_{\mm{N}}}E_{\mm{N}}(\tilde{\rho},\tilde{\mathbf{v}})<\infty$.
Moreover, from Proposition \ref{pro:nnn0203} we see that there is a pair of functions $({\rho},{\mathbf{v}})\in\mathcal{A}_{\mm{N}}$,
such that \begin{equation*}
E_{\mm{N}}({\rho},{\mathbf{v}})>0.\end{equation*}
Hence, $\sup_{(\tilde{\rho},\tilde{\mathbf{v}})\in\mathcal{A}_{\mm{N}}}E_{\mm{N}}
(\tilde{\rho},\tilde{\mathbf{v}})>0$.
Next we show that a maximizer of (\ref{n0206}) exists and that the corresponding
Euler-Lagrange equations are equivalent to the boundary problem
(\ref{n0201}).

\begin{pro}\label{npro:0201}
Assume that the steady density profile $\bar{\rho}$ satisfies \eqref{0102} and \eqref{nnnn0103}, then
there holds that
 \begin{enumerate}
    \item[(1)] $E_{\mm{N}}(\tilde{\rho},{\tilde{\mathbf{v}}})$
achieves its supremum on $\mathcal{A}_{\mm{N}}$.
   \item[(2)] Let $(\tilde{\rho}_0,\bar{\mathbf{v}}_0)$ be
a maximizer and $\Lambda=\sup_{(\tilde{\rho},\tilde{\mathbf{v}})\in
\mathcal{A}_{\mm{N}}}E_{\mm{N}}(\tilde{\rho},\tilde{\mathbf{v}})$, then there is a $\tilde{p}_0$
 such that ($\tilde{\rho}_0$, $\tilde{\mathbf{v}}_0$, $\tilde{p}_0$, $\Lambda$) satisfies
the boundary problem (\ref{n0201}). Moreover,
$(\tilde{\rho}_0,\tilde{\mathbf{v}}_0,\tilde{p}_0) \in H^1\times H^2\times H^1$.
 \end{enumerate}
\end{pro}
\begin{pf}
(1) 
Let $(\tilde{\rho}_n,\tilde{\mathbf{v}}_n)\in \mathcal{A}_{\mm{N}}$ be a maximizing sequence, then
$E_{\mm{N}}(\tilde{\rho}_n,\tilde{\mathbf{v}}_n)$ is lower bounded. This together with
(\ref{n0205}) implies that
$(\tilde{\rho}_n,\tilde{\mathbf{v}}_n)$ is bounded in $L^2\times H^\sigma_0$. Thus there exists a
$(\tilde{\rho}_0,\tilde{\mathbf{v}}_0)\in L^2\times H^\sigma_0$ and a subsequence (still denoted by
$(\tilde{\rho}_n,\tilde{\mf{v}}_n)$ for simplicity), such that $\tilde{\rho}_n\rightarrow \tilde{\rho}_0$ weakly in $L^2$,
  and $\tilde{\mathbf{v}}_n\rightarrow \tilde{\mathbf{v}}_0$ weakly in
$H^\sigma_0$ and strongly in $L^2$. Moreover, by the lower semi-continuity, one has
\begin{equation*}
\begin{aligned}
0<&\sup_{(\tilde{\rho},\tilde{\mathbf{v}})\in
\mathcal{A}_{\mm{N}}}E_{\mm{N}}(\tilde{\rho},\tilde{\mathbf{v}})
= \limsup_{n\rightarrow\infty}E_{\mm{N}}(\tilde{\rho}_n,\tilde{\mathbf{v}}_n)\\
= &- \frac{\mu}{g}\liminf_{n\rightarrow \infty}
\int|\nabla\tilde{\mf{v}}_n|^2\mm{d}\mf{x}-2\lim_{n\rightarrow
\infty}\int \tilde{\rho}_n\tilde{v}_{n3}\mm{d}\mf{x}\\
\leq & E(\tilde{\rho}_0,\tilde{\mathbf{v}}_0)\leq \sup_{(\tilde{\rho},\tilde{\mathbf{v}})\in
\mathcal{A}_{\mm{N}}}E_{\mm{N}}(\tilde{\rho},\tilde{\mathbf{v}}),
\end{aligned}\end{equation*}
and
\begin{equation*}J((\tilde{\rho}_0,\tilde{\mathbf{v}}_0))\leq 1,
\end{equation*}
which imply
\begin{equation*}
\begin{aligned}
\sup_{(\tilde{\rho},\tilde{\mathbf{v}})\in
\mathcal{A}_{\mm{N}}}E_{\mm{N}}(\tilde{\rho},\tilde{\mathbf{v}})
= E(\tilde{\rho}_0,\tilde{\mathbf{v}}_0)\;\mbox{ and }\; J((\tilde{\rho}_0,\tilde{\mathbf{v}}_0))>0,
\end{aligned}\end{equation*}
where we have denoted the third component of $\tilde{\mf{v}}_n$ by $\tilde{{v}}_{n3}$.

Suppose by contradiction that $J((\tilde{\rho}_0,\tilde{\mathbf{v}}_0))<1$. In view of the homogeneity
of $J$, we may find an $\alpha>1$ so that $J(\alpha (\tilde{\rho}_0,\tilde{\mathbf{v}}_0))=1$,
i.e., we may scale up
$(\tilde{\rho}_0,\tilde{\mathbf{v}}_0)$ so that $\alpha (\tilde{\rho}_0,\tilde{\mathbf{v}}_0)\in
\mathcal{A}_{\mm{N}}$. From this we deduce that
\begin{equation*}
E(\alpha(\tilde{\rho}_0,\tilde{\mathbf{v}}_0))=\alpha^2E(\tilde{\rho}_0,\tilde{\mathbf{v}}_0) =
\alpha^2\sup_{(\tilde{\rho},\tilde{\mathbf{v}})\in
\mathcal{A}_{\mm{N}}}E(\tilde{\rho},\tilde{\mathbf{v}})> \sup_{(\tilde{\rho},\tilde{\mathbf{v}})\in
\mathcal{A}_{\mm{N}}}E(\tilde{\rho},\tilde{\mathbf{v}}),
\end{equation*}
which is a contradiction since $\alpha (\tilde{\rho},\tilde{\mathbf{v}}_0)\in
\mathcal{A}_{\mm{N}}$. Hence $J((\tilde{\rho}_0,\tilde{\mathbf{v}}_0))=1$.
 This shows that $E(\tilde{\rho},\tilde{\mathbf{v}})$
achieves its  supremum on $\mathcal{A}_{\mm{N}}$.

(2) Notice that since $E(\tilde{\rho},\tilde{\mathbf{v}})$ and $J(\tilde{\rho},\tilde{\mathbf{v}})$
 are homogeneous of degree $2$, the definition (\ref{n0206}) is equivalent to
  \begin{equation}\label{n0227}\Lambda=\sup_{(\tilde{\rho},\tilde{\mathbf{v}})\in
  L^2\times {H}^{\sigma}_0}\frac{E(\tilde{\rho},\tilde{\mathbf{v}})}{J(\tilde{\rho},\tilde{\mathbf{v}})}.
\end{equation}
For any $\tau\in \mathbb{R}$, $\phi\in L^2$ and $\mf{w}\in {H}^\sigma_0$, we take
$(\tilde{\phi},\tilde{\mf{w}}):=(\tilde{\rho}_0+\tau\phi,
\tilde{\mathbf{v}}_0+\tau\mf{w})$, then (\ref{n0227}) implies
  \begin{equation*}E(\tilde{\phi},\tilde{\mf{w}})-\Lambda J(\tilde{\phi},\tilde{\mf{w}})\leq 0.
\end{equation*}
If we set
$I(\tau)=E(\tilde{\phi},\tilde{\mf{w}})-\Lambda J(\tilde{\phi},\tilde{\mf{w}})$,
then   $I(\tau)\in C^1(\mathbb{R})$, $I(\tau)\leq 0$ for all $\tau\in \mathbb{R}$ and
$I(0)=0$. This implies $I'(0)=0$. Hence, a direct computation leads to
  \begin{equation*}\begin{aligned}
\Lambda\int
     \left( \frac{\tilde{\rho}_0\phi }{   \bar{\rho}'} +     \frac{\bar{\rho}\tilde{\mathbf{v}}_0\cdot\mf{w}}{g}\right)
     \mm{d}\mf{x}=-\int \left(\frac{\mu}{g}\nabla\tilde{\mf{v}}_0:\nabla\mf{w}+
\tilde{\rho}_0{{w}}_3+\tilde{{v}}_{03}\phi\right)\mm{d}\mf{x},\end{aligned}  \end{equation*}
which implies that
  \begin{eqnarray*}
&&\Lambda\int
   \frac{\tilde{\rho}_0\phi                      }{   \bar{\rho}'} \mm{d}\mf{x}=-\int
\tilde{{v}}_{03}\phi\mm{d}\mf{x}\quad\mbox{ for any }\phi\in L^2,\\
&&\Lambda\int
     \frac{\bar{\rho}\tilde{\mathbf{v}}_0\cdot\mf{w}}{g}\mm{d}\mf{x}=-\int \left( \frac{\mu}{g}\nabla\tilde{\mf{v}}_0:\nabla\mf{w}+
\tilde{\rho}_0\mf{e}_3\cdot{\mf{w}}\right)\mm{d}\mf{x}\quad\mbox{ for any }\mf{w}\in {H}^\sigma_0,
\end{eqnarray*}
i.e.,
$(\tilde{\rho}_0,\tilde{\mathbf{v}}_0)$ is a weak
solution to the boundary problem \eqref{n0201}. In particular, $\Lambda\tilde{\rho}_0/\bar{\rho}'=-\tilde{{v}}_{03}$
and
  \begin{eqnarray*}
\int
     \nabla\tilde{\mf{v}}_0:\nabla\mf{w}\mm{d}\mf{x}=\frac{1}{{\mu}\Lambda}\int \left( {g\bar{\rho}'\tilde{v}_{03}}\mf{e}_3-
    {\Lambda^2\bar{\rho}\tilde{\mathbf{v}}_0}
\right)\cdot\mf{w}\mm{d}\mf{x}\quad\mbox{ for any }\mf{w}\in {H}^\sigma_0.
\end{eqnarray*}
 Thus, it follows from the classical regularity theory on the Stokes equations that there are constants $c_6$ dependent of
 the domain $\Omega$, and $c_7$ dependent of $g$, $\Lambda$, $\Omega$ and $\bar{\rho}$, such that
  \begin{equation*}  
  \|\tilde{\mathbf{v}}\|_{H^{2}}+   \|\nabla \tilde{p}_0\|_{L^{2}}\leq c_6\|
(g\bar{\rho}'\tilde{{v}}_{03}\mf{e}_3-\Lambda^2\bar{\rho}\tilde{\mathbf{v}}_0)/(\mu\Lambda )\|_{L^2}\leq c_7,
\end{equation*}
where $\tilde{p}_0\in H^{1}$ is the corresponding pressure field associated to $\tilde{\mathbf{v}}_0$. Hence,
$(\tilde{\rho}_0,\tilde{\mathbf{v}}_0,\tilde{p}_0) \in H^1\times H^2\times H^1$ is a strong
solution to the boundary problem \eqref{n0201} with $\Lambda:=\sup_{(\tilde{\rho},\tilde{\mathbf{v}})\in
\mathcal{A}_{\mm{N}}}E_\mm{N}(\tilde{\rho},\tilde{\mathbf{v}})$. This completes the proof.
\hfill $\Box$
\end{pf}

Now we are able to show that the growth rate constructed by \eqref{growth} is indeed equal to that
defined by \eqref{n0206}.
\begin{pro}\label{equal}
Let $\Lambda$ be constructed by \eqref{growth}, and \begin{equation*}
\Lambda_\mm{N}:=\sup_{(\tilde{\rho},\tilde{\mathbf{v}})\in
\mathcal{A}_{\mm{N}}}E_\mm{N}(\tilde{\rho},\tilde{\mathbf{v}}),\end{equation*}
then $\Lambda=\Lambda_N$.
\end{pro}
\begin{pf}
In view of Proposition \ref{pro:nnn0203}, there are functions ($\tilde{\rho}_0$, $\tilde{\mathbf{v}}_0$, $\tilde{p}_0)\in H^1\times H^2\times H^1$
which solve the boundary problem \eqref{n0201} with $\Lambda$ constructed by \eqref{growth}. Hence,
\begin{equation*} \begin{aligned}
                            \Lambda = -{
\int\left(\frac{\mu |\nabla\tilde{\mf{v}}_0|^2}{g}+2\tilde{\rho}_0\tilde{v}_{03}\right)\mm{d}\mf{x}}
\bigg/{\int
\left( \frac{\tilde{\rho}^2_0}{\bar{\rho}'}+\frac{\bar{\rho}|\tilde{\mathbf{v}}_0|^2}{g}\right)
\mm{d}\mf{x}}.
\end{aligned}\end{equation*}
By virtue of the dentition \eqref{n0206}, $\Lambda\leq \Lambda_N$.

On the other hand, by Proposition \ref{npro:0201}, we see that there are functions ($\bar{\rho}_0$, $\bar{\mathbf{v}}_0$, $\bar{p}_0$)
which solve the boundary problem \eqref{n0201} with $\Lambda=\Lambda_{\mm{N}}>0$. Hence,
 \begin{equation*}\begin{aligned}
\label{0111}0< \Lambda_{\mm{N}}^2
=\frac{g \int_{\Omega}\bar{\rho}' \bar{
{v}}_3^2\mm{d}\mf{x}-\Lambda_{\mm{N}}\mu
\int_{\Omega}|\nabla
\bar{\mathbf{v}}|^2\mm{d}\mf{x}}
{\int_{\Omega}\bar{\rho}|\bar{\mathbf{v}}|^2\mathrm{d}\mf{x}}\leq \sup_{\tilde{\mathbf{v}}\in
\mathcal{A}}E(\tilde{\mathbf{v}},\Lambda_{\mm{N}} ).
\end{aligned}\end{equation*}
Thanks to \eqref{0240}, $\Lambda_\mm{N}\in (0,\mathfrak{S})$. Thus, we can use the monotonicity of
$\sup_{\tilde{\mathbf{v}}\in\mathcal{A}}E(\tilde{\mathbf{v}},s)$ on the variable $s$ to infer that
 \begin{equation*}\Lambda_{\mm{N}}^2\leq
\sup_{\tilde{\mathbf{v}}\in
\mathcal{A}}E(\tilde{\mathbf{v}},\Lambda_{\mm{N}} )
\leq \sup_{\tilde{\mathbf{v}}\in
\mathcal{A}}E(\tilde{\mathbf{v}},\Lambda )=\Lambda^2.
\end{equation*}
Hence, $\Lambda_{\mm{N}}\leq \Lambda $, which gives $\Lambda=\Lambda_{\mm{N}}$.  \hfill $\Box$
\end{pf}

Next, we use the energy method to show that $\Lambda$ is the sharp growth rate of arbitrary solutions
to the linearized problem (\ref{0108}).
\begin{pro}\label{grwothe} Assume that the assumption of Theorem \ref{thm:0101} is satisfied, and
let $(\varrho,\mf{u},q)$ solve the following linearized problem
\begin{equation}\label{0301}\left\{
  \begin{array}{l}
 \varrho_t+ \bar{\rho}'{u}_3=0, \\
    \bar{\rho}\mathbf{u}_t +\nabla q+g\varrho e_3=\mu \Delta \mathbf{u},\\
\mm{div}\mf{u}=0
  \end{array}
\right.
\end{equation}
with initial and boundary conditions
$(\varrho,{\bf u})|_{t=0}=(\varrho_0,{\bf u}_0)$ in $\Omega$ and $
{\bf u}|_{\partial\Omega}={\bf 0}$ for any $t>0$.
Then, we have the following estimates for any $t\geq 0$.
\begin{eqnarray}
&&\label{0313} \|\varrho(t)\|_{X}^2\leq Ce^{2\Lambda
t}(\|\varrho_0\|_{X}^2+\|\mathbf{u}_0\|_{H^2}^2),\quad
X=L^2\mbox{  or }H^1,\\[1mm]
&&\label{uestimate} \|\mathbf{u}(t)\|_{H^2 }^2+\|\mathbf{u}_t(t)\|^2_{L^2 }+
\int_0^t\|\nabla\mathbf{u}(s)\|^2_{L^2}\mm{d}s\leq
Ce^{2\Lambda t}(\|\varrho_0\|_{L^2}^2+\|\mathbf{u}_0\|_{H^2}^2),
\end{eqnarray}
where $\Lambda$ is constructed by \eqref{growth}, and the constant $C$ may depend on $g$, $\mu$, $\bar{\rho}$ and $\Lambda$. In addition,
if $\bar{\rho}$ further satisfies $\inf_{\mf{x}\in \Omega}\{\bar{\rho}'(\mf{x})\} >0$, then
\begin{equation}\label{L2es}
\|(\varrho,\mathbf{u})(t)\|_{L^2}\leq Ce^{\Lambda t}\|(\varrho_0,\mathbf{u}_0)\|_{L^2},\end{equation}
where the constant $C$ only depends on $g$ and $\bar{\rho}$.
\end{pro}
\begin{pf}
Let
 $(\varrho,\mathbf{u},q)$ solve the initial and boundary value problem above.
 In what follows, we denote by $C$ a generic positive constant which may depend on $\mu$, $g$, $\Lambda$, $\bar{\rho}$
and $\Omega$.

We differentiate the linearized  momentum (\ref{0301})$_{2}$ in time, multiply the resulting
equation by $\mf{u}_t$ and then integrate (by parts) over $\Omega$ to obtain
\begin{equation*} \frac{d}{dt}
\int\bar{\rho}|\mathbf{u}_t|^2
\mathrm{d}\mathbf{x}+2 g\int{\varrho}_t\partial_t{u}_3
\mathrm{d}\mathbf{x}+
2\mu\int|\nabla\partial_t\mathbf{u}|^2\mm{d}\mf{x}=0.
\end{equation*}
Using the linearized mass equation (\ref{0301})$_1$, we get
\begin{equation}\label{nnn0314}\frac{d}{dt}
\int \left(\bar{\rho}|\mathbf{u}_t|^2
-g\bar{\rho}'{u}_3^2\right)\mathrm{d}\mathbf{x}+
2\mu\int |\nabla\partial_t\mathbf{u}|^2\mm{d}\mf{x}=0.
\end{equation}
Thanks to \eqref{growth}, we have
\begin{equation}\label{0302}\int  g\bar{\rho}'|u_3|^2\mm{d}\mf{x}-
{\Lambda}\mu\int
|\nabla\mathbf{u}|^2\mm{d}\mf{x}\leq {\Lambda^2}\int \bar{\rho}|\mf{u}|^2\mm{d}\mf{x}.
\end{equation}
Thus, integrating (\ref{0301}) in time from $0$ to $t$ and using \eqref{0302}, we get
\begin{equation}\label{0314}
\|\sqrt{\bar{\rho}} \mathbf{u}_t(t)\|^2_{L^2}+2\mu\int_0^t\|\nabla\partial_s\mathbf{u}(s)\|^2_{L^2}\mm{d}s
\leq I_0+{\Lambda^2} \|\sqrt{\bar{\rho}}\mf{u}(t)\|_{L^2}+ {\Lambda}\mu\|\nabla\mathbf{u}(t)\|^2_{L^2}
\end{equation}
with
\begin{equation*}
I_0=\|\sqrt{\bar{\rho}} \mathbf{u}_t|_{t=0}\|_{L^2}^2-
g\|u_3(0)\|^2_{L^2},\end{equation*}
where the first term on the right-hand side of (\ref{0314}) can bounded in terms of
the initial data $(\varrho_0,\mf{u}_0)$. In fact, from \eqref{0301}$_2$ and \eqref{0102} we get
   \begin{equation*}  \begin{aligned}
    \int \bar{\rho}|{\bf u}_\tau|^2(\tau)
    \mm{d}\mf{x} = & \int [\mu \Delta \mathbf{u}-g\varrho \mf{e}_3]\cdot\mf{u}_{\tau}(\tau)\mm{d}\mf{x}\\
\leq & C(\|\varrho (\tau)\|_{L^2}^2+\|\Delta \mf{u}(\tau)\|_{L^2}^2 +\frac{\inf{\bar{\rho}}}{2} \|{\bf u}_\tau(\tau)\|^2_{L^2},
\end{aligned} \end{equation*}
which implies
 \begin{equation*} \begin{aligned}
      \| {\bf u}_\tau(\tau)\|^2_{L^2}
\leq C(\|\varrho(\tau)\|_{L^2}^2+\|\Delta \mf{u}(\tau)\|_{L^2}^2
\;\;\mbox{ for any }\tau>0.  \end{aligned} \end{equation*}
Let $\tau\rightarrow 0$, we see that
 \begin{equation*}\label{initial}\| {\bf u}_t|_{t=0}\|^2_{L^2}
 \leq  C(\|\varrho_0\|_{L^2}^2+\|\Delta\mf{u}_0\|_{L^2}^2).
 \end{equation*}

Integrating in time and using Cauchy-Schwarz's inequality, we find that
 \begin{equation}\begin{aligned}\label{0316}
\Lambda\|\nabla\mathbf{u}(t)\|^2_{L^2}=&
\Lambda\|\nabla\mathbf{u}_0\|^2_{L^2}+2\Lambda\int_0^t
\int_{\Omega}\nabla\partial_s\mathbf{u}(s):\nabla\mathbf{u}(s)\mm{d}\mf{x}\mathrm{d}s\\
\leq &\Lambda\|\nabla\mathbf{u}_0\|^2_{L^2}+\int_0^t\|\nabla\partial_s\mathbf{u}(s)\|^2_{L^2}\mathrm{d}s
+\Lambda^2\int_0^t\|\nabla\mathbf{u}(s)\|^2_{L^2}\mathrm{d}s.
\end{aligned}\end{equation}
On the other hand,
\begin{equation}\begin{aligned}\label{0317}
\Lambda\partial_t\|\sqrt{\bar{\rho}}\mathbf{u}(t)\|^2_{L^2}=2\Lambda\int_{\Omega}
\bar{\rho}\mathbf{u}(t)\cdot\partial_t\mathbf{u}(t)\mm{d}\mf{x}\leq\|\sqrt{\bar{\rho}}\partial_t\mathbf{u}(t)\|^2_{L^2}
+\Lambda^2\|\sqrt{\bar{\rho}}\mathbf{u}(t)\|^2_{L^2}.
\end{aligned}\end{equation}
Hence, putting (\ref{0314})--(\ref{0317}) together, we obtain the differential inequality
$$ \partial_t\|\sqrt{\bar{\rho}}\mathbf{u}(t)\|^2_{L^2}+\mu\|\nabla\mathbf{u}(t)\|^2_{L^2}
\leq I_1+2\Lambda\left(\|\sqrt{\bar{\rho}}\mathbf{u}(t)\|^2_{L^2}+\mu\int_0^t\|\nabla\mathbf{u}(s)\|^2_{L^2}
\mathrm{d}s\right)$$
with $I_1=2\mu\|\nabla\mathbf{u}_0\|_{L^2}+I_0/\Lambda$. An application of Gronwall's inequality
then shows that
\begin{equation}\begin{aligned}\label{0319}
\|\sqrt{\bar{\rho}}\mathbf{u}(t)\|^2_{L^2}+\mu\int_0^t\|\nabla\mathbf{u}(s)\|^2_{L^2}\mm{d}s
\leq & e^{2\Lambda t}\|\sqrt{\bar{\rho}}\mathbf{u}_0\|^2_{L^2}
+\frac{I_1}{2\Lambda}\left(e^{2\Lambda t}-1\right)\\
\leq & Ce^{2\Lambda t}(\|\varrho_0\|_{L^2}^2+
\|\mathbf{u}_0\|_{H^2}^2)
\end{aligned}\end{equation}
for any $t\geq 0$. Thus, making use of (\ref{0314}), (\ref{0316}) and (\ref{0319}), we deduce that
\begin{equation*}\begin{aligned}
\frac{1}{\Lambda}\|\sqrt{\bar{\rho}} \mathbf{u}_t(t)\|^2_{L^2}+
\mu\|\nabla\mathbf{u}(t)\|^2_{L^2} & \leq  I_1+ {\Lambda}\|\sqrt{\bar{\rho}}\mathbf{u}(t)\|^2_{L^2}+2
{\Lambda}\mu\int_0^t\|\nabla\mathbf{u}(s)\|^2_{L^2}\mm{d}s \\
&\leq Ce^{2\Lambda t}(\|\varrho_0\|_{L^2}^2+\|\mathbf{u}_0\|_{H^2}^2),
\end{aligned}\end{equation*}
which, together with Poincar\'e's inequality, yields
\begin{eqnarray}\label{uestimate1}
\|\mathbf{u}(t)\|_{H^1 }^2+\|\mathbf{u}_t(t)\|^2_{L^2 }+
\int_0^t\|\nabla\mathbf{u}(s)\|^2_{L^2}\mm{d}s\leq
Ce^{2\Lambda t}(\|\varrho_0\|_{L^2}^2+\|\mathbf{u}_0\|_{H^2}^2),
\end{eqnarray}
Hence, using (\ref{0108})$_1$, we have
\begin{equation}\begin{aligned}\label{0323}
\|\varrho(t)\|_{X}\leq&
\|\varrho_0\|_{X}+\int_0^t\|\varrho_s(s)\|_{X}\mm{d}s\\
\leq &\|\varrho_0\|_{X}+\|\bar{\rho}'\|_{L^\infty}\int_0^t\|u_3(s)\|_{X}\mm{d}s\\
\leq & Ce^{\Lambda t}(\|\varrho_0\|_{X}+\|\mathbf{u}_0\|_{H^2}),
\end{aligned}\end{equation}
where $X=L^2$ or $H^1$. Thus (\ref{0313}) is proved.

From the classical regularity theory on the Stokes equations it follows that
\begin{equation*}
\begin{aligned}\|\nabla^2 \mathbf{u} \|_{L^2}^2+\|\nabla q\|_{L^2}^2
\leq C\|\bar{\rho} \mathbf{u}_t+{g}\varrho \mf{e}_3\|_{L^2}^2
\leq  C (\|\mathbf{u}_t\|_{L^2}^2+\|\varrho\|_{L^2}^2),
\end{aligned}\end{equation*}
which combined with \eqref{uestimate1} and \eqref{0313} gives \eqref{uestimate}.

Finally, we prove \eqref{L2es} under the condition of $\inf_{\mf{x}\in \Omega}\{\bar{\rho}'(\mf{x})\}>0$.
Multiplying \eqref{0301}$_1$ and \eqref{0301}$_2$ by $\varrho/\bar{\rho}'$ and $\mf{u}/g$, respectively,
adding the resulting equalities, and integrating over $\Omega$ (by parts), one has
 \begin{equation*}
\begin{aligned}
 \frac{1}{2}\frac{d}{dt}\int\left(\frac{|\varrho(t)|^2}{\bar{\rho}'} + \frac{\bar{\rho}|\mathbf{u}(t)|^2}{g}\right)\mm{d}\mf{x}
 +\int \left(\frac{\mu|\nabla \mathbf{u}|^2}{g}+2\varrho u_3 \right)\mm{d}\mf{x} =0.
 \end{aligned}  \end{equation*}
On the other hand, thanks to \eqref{n0206} and Proposition \ref{equal}, we have
\begin{equation}\label{supress}
\begin{aligned}-\int \left(\frac{\mu|\nabla \mathbf{u}
|^2}{g}+2\varrho u_3 \right)\mm{d}\mf{x}\leq
\Lambda\int\left(\frac{|\varrho|^2}{\bar{\rho}'} + \frac{\bar{\rho}|\mathbf{u}|^2}{g}\right)\mm{d}\mf{x}.
 \end{aligned}  \end{equation}
Consequently, one has
   \begin{equation*}
\begin{aligned}
\frac{d}{dt}\int\left(\frac{|\varrho(t)|^2}{\bar{\rho}'} +  \frac{\bar{\rho}|\mathbf{u}(t)|^2}{g}\right)\mm{d}\mf{x}
\leq 2\Lambda\int\left(\frac{|\varrho|^2}{\bar{\rho}'} + \frac{\bar{\rho}|\mathbf{u}|^2}{g}\right), \end{aligned}  \end{equation*}
which together with Gronwall's inequality yields
 \begin{equation*}
\begin{aligned}
 \int\left(\frac{|\varrho(t) |^2}{\bar{\rho}'} +
 \frac{\bar{\rho}|\mathbf{u}(t)|^2}{g}\right)\mm{d}\mf{x}
&\leq e^{2\Lambda t}\int\left(\frac{|\varrho_0 |^2}{\bar{\rho}'} +
 \frac{\bar{\rho}|\mathbf{u}_0|^2}{g}\right)\mm{d}\mf{x}.
 \end{aligned}  \end{equation*}
 Therefore, \eqref{L2es} is obtained.    \hfill $\Box$
\end{pf}

\section{Nonlinear energy estimates}\label{sec:04}

In this section, we first formally derive some  nonlinear energy
estimates for the perturbed problem, which will  be used in the
proof of Theorem \ref{thm:0102} in Section \ref{sec:05}. To this end, let
$(\varrho,\mathbf{u},q)$ be a classical solution of the perturbed
problem (\ref{0105})--(\ref{0107}) with $\rho:=\varrho+\bar{\rho}>0$,
and the initial density satisfies
$0<\inf_{\mathbf{x}\in\Omega}\{\rho_0(\mathbf{x})\}\leq \sup_{\mathbf{x}\in\Omega}\{\rho_0({\mathbf{x}})\}<\infty$,
where $\rho_0:=\varrho_0+\bar{\rho}$. Moreover, we assume
that the classical solution $(\varrho,\mathbf{u},q)$ is sufficiently regular so that
the procedure of formal deduction makes sense. In what follows, we denote by $C$ a
generic positive constant which may depend on $\mu$, $g$, $\bar{\rho}$ and $\Omega$.

\subsection{Estimates for the density and velocity}

 We first observe that the continuity equation (\ref{0105})$_1$ and the incompressibility condition
(\ref{0105})$_3$ imply immediately that for any $t>0$ and $\mf{x}\in \Omega$,
\begin{equation}\label{0402}0<\inf_{\mathbf{x}\in\Omega}\{\rho_0(\mathbf{x})\}
\leq \rho(t,\mf{x})\leq \sup_{\mathbf{x}\in\Omega}\{\rho_0({\mathbf{x}})\}.
\end{equation}
This property will be repeatedly used to bound the density below.

From the continuity equation and the incompressibility condition it follows that
\begin{equation}\label{0403}\frac{d}{dt}
\|\varrho(t)\|_{L^2}^2=-2\int \bar{\rho}'\varrho u_3\mathrm{d}\mathbf{x}.
\end{equation}
Multiplying (\ref{0105})$_2$ by $\mf{u}$, using (\ref{0105})$_1$,
and then integrating (by parts) over $\Omega$, we obtain
\begin{equation*}
\frac{d}{dt}\int\rho|\mathbf{u}|^2(t)\mathrm{d}\mathbf{x}+2\mu\int|\nabla \mathbf{u}|^2
\mathrm{d}\mathbf{x}= -2g\int {\varrho}{u}_3\mathrm{d}\mathbf{x}.
\end{equation*}
Adding the above two equalities and using H\"older's inequality, we conclude
\begin{equation}\label{0406}\begin{aligned}
\frac{d}{dt}\|(\varrho,\sqrt{{\rho}}\mathbf{u})(t)\|_{L^2}^2+2\mu\|\nabla \mathbf{u}\|_{L^2}^2
\leq & C\left\|\varrho\right\|_{L^2}\|\mathbf{u}\|_{L^2}.
\end{aligned}\end{equation}

To control $\mathbf{u}_t$, we multiply (\ref{0105})$_2$ by $\mathbf{u}_t$ in $L^2$ and
apply Cauchy-Schwarz's inequality  to infer that
\begin{equation}\label{nnn0405} \begin{aligned}\|\sqrt{\rho}\mathbf{u}_t\|_{L^2}^2+\mu\frac{d}{dt}
\|\nabla \mathbf{u}(t)\|_{L^2}^2\leq &C\|({ \varrho},\sqrt{\rho}\mathbf{u}\cdot\nabla\mathbf{u})\|_{L^2}^2,
\end{aligned}\end{equation}
where the second term on the right hand side can be bounded as follows, using
the embedding theorem, and Poincar$\mm{\acute{e}}$'s, H\"older's and interpolation inequalities.
\begin{equation*}
\|\mathbf{u}\cdot\nabla
\mathbf{u}\|_{L^2}^2\leq \|\mathbf{u}\|_{L^6}^2\|\nabla
\mathbf{u}\|_{L^3}^2\leq \|\mathbf{u}\|_{L^6}^2\|\nabla
\mathbf{u}\|_{L^2}\|\nabla
\mathbf{u}\|_{L^6}\leq C\|\nabla
\mathbf{u}\|_{L^2}^3\|\nabla
\mathbf{u}\|_{H^1}.\end{equation*}
Putting the above inequality into \eqref{nnn0405} and using \eqref{0402}, we conclude that
\begin{equation}\label{0410} \begin{aligned}\|\sqrt{\rho}\mathbf{u}_t\|_{L^2}^2+\mu\frac{d}{dt}
\|\nabla \mathbf{u}(t)\|_{L^2}^2 \leq C(\|{ \varrho}\|_{L^2}^2 +\|\nabla \mathbf{u}\|_{L^2}^3
\|\nabla\mathbf{u}\|_{H^1}).\end{aligned}\end{equation}

We proceed to derive higher-order estimates of the density and velocity.
Using (\ref{0105})$_1$ and keeping in mind that $p=q+\bar{p}$, we can rewrite (\ref{0105})$_2$ as
\begin{equation*}
\rho\mathbf{u}_t+\rho\mathbf{u}\cdot \nabla \mathbf{u}+\nabla
p=\mu\Delta\mathbf{u}-\rho g e_3,
\end{equation*}
whence, by taking the time derivative,
\begin{equation*}
\rho\mathbf{u}_{tt}+\rho\mathbf{u}\cdot \nabla
\mathbf{u}_t-\mu\Delta\mathbf{u}_t+\nabla
p_t=-\rho_t(\mathbf{u}_t+\mathbf{u}\cdot\nabla \mathbf{u}+g
e_3)-\rho\mathbf{u}_t\cdot\nabla \mathbf{u},
\end{equation*}
which, by using the continuity equation, can be written as
\begin{equation*} \begin{aligned}
&\rho\left(\frac{1}{2}|\mathbf{u}_{t}|^2\right)_t+\rho\mathbf{u}\cdot
\nabla
\left(\frac{1}{2}|\mathbf{u}_t|^2\right)-\mu\Delta\mathbf{u}_t\cdot\mathbf{u}_t+\nabla
p_t\cdot\mathbf{u}_t\\
&=\mathrm{div}(\rho\mathbf{u})(\mathbf{u}_t+\mathbf{u}\cdot\nabla
\mathbf{u}+ g e_3)\cdot\mathbf{u}_t-\rho(\mathbf{u}_t\cdot\nabla
\mathbf{u})\mathbf{u}_t.\end{aligned}
\end{equation*}
Hence, by integrating by parts, we see that
\begin{equation}\label{0425}\begin{aligned}
&\frac{1}{2}\frac{d}{dt}\int\rho|\mathbf{u}_{t}(t)|^2\mathrm{d}\mathbf{x}
+\mu\int|\nabla \mathbf{u}_t|^2\mathrm{d}\mathbf{x}\\
&\leq \int2\rho|\mathbf{u}||\mathbf{u}_t||\nabla
\mathbf{u}_t|+\rho|\mathbf{u}||\mathbf{u}_t||\nabla
\mathbf{u}|^2+\rho |\mathbf{u}|^2|\mathbf{u}_t||\nabla^2 \mathbf{u}|\\
&\quad +\rho|\mathbf{u} |^2|\nabla \mathbf{u}||\nabla
\mathbf{u}_t|+\rho|\mathbf{u}_t|^2|\nabla
\mathbf{u}|+g\rho|\mathbf{u}||\nabla
\mathbf{u}_t|:=\sum_{i=1}^6I_j,\end{aligned}
\end{equation}
where $I_j$ can be bounded as follows, employing straightforward calculations.
\begin{eqnarray*}
&& I_1\leq
C\|\mathbf{u}\|_{L^6}\|\mathbf{u}_t\|_{L^3}\|\nabla
\mathbf{u}_t\|_{L^2}\leq C\|\mathbf{u}\|_{L^6}\|\mathbf{u}_t\|_{L^2}^{1/2}
\|\mathbf{u}_t\|_{L^6}^{1/2}\|\nabla\mathbf{u}_t\|_{L^2}\\
&& \quad \leq C\|\nabla
\mathbf{u}\|_{L^2}\|\mathbf{u}_t\|_{L^2}^{1/2}\|\nabla
\mathbf{u}_t\|_{L^2}^{3/2}\leq C(\varepsilon)\|\nabla
\mathbf{u}\|_{L^2}^4\|\mathbf{u}_t\|_{L^2}^2+\varepsilon\|\nabla\mathbf{u}_t\|_{L^2}^2,\\[1mm]
&& I_2\leq C\|\nabla \mathbf{u} \|_{L^2}\|\nabla
\mathbf{u}_t\|_{L^2}\|\nabla \mathbf{u}\|_{L^2}\|\nabla
\mathbf{u}\|_{H^1}\leq C(\varepsilon)\|\nabla
\mathbf{u}\|_{L^2}^4\|\nabla
\mathbf{u}\|_{H^1}^2+\varepsilon\|\nabla\mathbf{u}_t \|_{L^2}^2,\\[1mm]
&&I_3\leq C\|\nabla \mathbf{u}\|_{L^2}^2\|\nabla
\mathbf{u}_t\|_{L^2}\|\nabla^2 \mathbf{u}\|_{L^2}\leq
C(\varepsilon)\|\nabla \mathbf{u}\|_{L^2}^4\|\nabla^2
\mathbf{u}\|_{L^2}^2+\varepsilon\|\nabla \mathbf{u}_t\|_{L^2}^2,\\[1mm]
&& I_4\leq C\|\nabla \mathbf{u}\|_{L^2}^2\|\nabla
\mathbf{u}\|_{H^1}\|\nabla \mathbf{u}_t\|_{L^2}\leq
C(\varepsilon)\|\nabla \mathbf{u}\|_{L^2}^4\|\nabla
\mathbf{u}\|_{H^1}^2+\varepsilon\|\nabla \mathbf{u}_t\|_{L^2}^2,\\
&& I_5\leq C\|\nabla
\mathbf{u}\|_{L^2}\|\mathbf{u}_t\|_{L^2}^{1/2}\|\nabla
\mathbf{u}_t\|_{L^2}^{3/2}\leq C(\varepsilon)\|\nabla
\mathbf{u}\|_{L^2}^4\|\mathbf{u}_t\|_{L^2}^2+\varepsilon\|\nabla
\mathbf{u}_t\|_{L^2}^2, \\[1mm]
&& I_6\leq C(\varepsilon)\|\mathbf{u}\|_{L^2}^2+\varepsilon\|\nabla\mathbf{u}_t\|_{L^2}^2.
\end{eqnarray*}
Here $C(\varepsilon)$ is a positive constant which depends on $\varepsilon$ and $\|\rho\|_{L^\infty}$.
Taking $\varepsilon=\mu/12$ and inserting all the above estimates into (\ref{0425}), we conclude
\begin{equation}\label{jfjs0427}
\frac{d}{dt}\|\sqrt{\rho}\mathbf{u}_{t}(t)\|_{L^2}^2+\mu\|\nabla
\mathbf{u}_t\|_{L^2}^2 \leq C\Big(\|\nabla
\mathbf{u}\|^4_{L^2}(\|\mathbf{u}_t\|^2_{L^2}+\|\nabla
\mathbf{u}\|_{H^1}^2)+\|\mathbf{u}\|_{L^2}^2\Big).
\end{equation}

\subsection{Energy estimates}
Adding \eqref{0406}, \eqref{0410} and \eqref{jfjs0427} together, we find that
\begin{equation*} \begin{aligned}
&\frac{d}{dt}
\|(\varrho,\sqrt{\rho}\mathbf{u},\sqrt{\mu}\nabla
\mathbf{u},\sqrt{\rho}\mathbf{u}_{t})(t)\|_{L^2}^2 +\|(\sqrt{2\mu}\nabla
\mathbf{u},\sqrt{\rho}\mathbf{u}_t,\sqrt{\mu}\nabla \mathbf{u}_t)\|_{L^2}^2 \\
&\leq C\left[\|(\varrho,\mathbf{u})\|_{L^2}^2
+\|\nabla \mathbf{u}\|^3_{L^2}(\|\nabla\mathbf{u}\|_{H^1}+\|\nabla
\mathbf{u}\|_{L^2}\|\mathbf{u}_t\|^2_{L^2}+\|\nabla\mathbf{u}\|_{L^2}\|\nabla\mathbf{u}\|_{H^1}^2)\right].
\end{aligned}\end{equation*}
Now if we denote
$${\mathcal{E}}(t)=\mathcal{E}((\varrho,\mathbf{u})(t))=\sqrt{
\|\varrho(t)\|_{L^2}^2+\|\mathbf{u}(t)\|_{H^2}^2}$$
and $$
\mathcal{E}_0=:\mathcal{E}(0)={{\mathcal{E}}(\varrho_0,\mathbf{u}_0)=\sqrt{
\|\varrho_0\|_{L^2}^2+\|\mathbf{u}_0\|_{H^2}^2}},$$
then there is a $\bar{\delta}_0\in (0,1]$ dependent of $\mu$, $\bar{\rho}$, $g$ and $\Omega$,
such that for any $\mathcal{E}(t)\leq \bar{\delta}_0$,
\begin{equation}\label{enegetinnnn}
\frac{d}{dt}\|(\varrho,\sqrt{\rho}\mathbf{u},\sqrt{\mu}\nabla\mathbf{u},\sqrt{\rho}\mathbf{u}_{t})(t)\|_{L^2}^2
 +\frac{1}{2}\|(\sqrt{2\mu}\nabla\mathbf{u},\sqrt{\rho}\mathbf{u}_t,\sqrt{\mu}\nabla
\mathbf{u}_t)\|_{L^2}^2  \leq C\|(\varrho,\mathbf{u})(t)\|_{L^2}^2.
\end{equation}

On the other hand, multiplying (\ref{0105})$_2$ by $\mathbf{u}_t$
in $L^2$ and recalling $\mathrm{div}\,\mathbf{u}_t=0$, we see that
\begin{equation*} \begin{aligned}
\int\rho|\mathbf{u}_t|^2(t)\mathrm{d}\mathbf{x}=\int(-\varrho{g}e_3-\rho\mathbf{u}\cdot \nabla
\mathbf{u}+\mu \Delta \mathbf{u})\cdot
\mathbf{u}_t\mathrm{d}\mathbf{x}, \end{aligned}
\end{equation*}
whence,
\begin{equation*}\label{0430}\begin{aligned}
\int\rho|\mathbf{u}_t|^2(t)\mathrm{d}\mathbf{x}\leq
C\int (\varrho^2+|\mathbf{u}|^2|\nabla
\mathbf{u}|^2+|\Delta\mathbf{u}|^2)(t)\mathrm{d}\mathbf{x}\leq C\mathcal{E}^2(t).\end{aligned}
\end{equation*}
Taking $t\to 0$ in the above inequality, one gets
\begin{equation*} \begin{aligned}
\limsup_{t\to 0}\int\rho|\mathbf{u}_t|^2(t)\mathrm{d}\mathbf{x}
\leq C\mathcal{E}_0^2.\end{aligned}
\end{equation*}
Therefore, we conclude from \eqref{enegetinnnn} that
\begin{equation}\label{enegetin123}
\|(\varrho,\mathbf{u},\nabla\mathbf{u},\mathbf{u}_{t})(t)\|_{L^2}^2
+ \int_0^t\|(\nabla\mathbf{u},\mathbf{u}_t,\nabla\mathbf{u}_t)(s)\|_{L^2}^2\mm{d}s
\leq C\left(\mathcal{E}_0^2+\int_0^t\|(\varrho,\mathbf{u})(s)\|_{L^2}^2\mm{d}s\right).
\end{equation}

We recall again that the pair ($\mathbf{u},q$) solves the Stokes equations:
\begin{equation*}
-\mu \Delta\mathbf{u} +\nabla q =-\rho \mathbf{u}_t-\rho (\mathbf{u}\cdot \nabla \mathbf{u})-\varrho
{g}e_3,\quad \mathrm{div}\mathbf{u}=0\;\; \mbox{ in }\Omega. \end{equation*}
Thus, it follows from the classical regularity theory on the Stokes equations that
\begin{equation}\label{0433}
\begin{aligned}\|(\nabla^2 \mathbf{u} ,\nabla q)\|_{L^2}^2
\leq &C\|-\rho \mathbf{u}_t-\rho(\mathbf{u}\cdot\nabla \mathbf{u})-\varrho {g}e_3\|_{L^2}^2\\
\leq & C( \|\mathbf{u}_t\|_{L^2}^2+ \|\mathbf{u}\|_{L^6}^2\|\nabla
\mathbf{u}\|_{L^3}^2+\|\varrho\|_{L^2}^2)\\
\leq & C\big[\|(\varrho,\mathbf{u}_t)\|_{L^2}^2+\|\nabla
\mathbf{u}\|_{L^2}^6\big] +\frac{1}{2}\|\nabla \mathbf{u}\|_{H^1}^2.
\end{aligned}\end{equation}
Hence, provided $\mathcal{E}(t)\leq \bar{\delta}_0\leq 1$,
\begin{equation*}
\begin{aligned}\frac{1}{2}\|\nabla^2 \mathbf{u} \|_{L^2}^2+\|\nabla q\|_{L^2}^2
\leq C\|(\varrho,\nabla
\mathbf{u},\mathbf{u}_t)\|_{L^2}^2,
\end{aligned}\end{equation*}
which, together with \eqref{enegetin123}, results in
\begin{equation}\label{enegetin1234}\begin{aligned}
& \mathcal{E}^2(t)+\|(\mathbf{u}_{t},\nabla q)(t)\|_{L^2}^2
+\int_0^t\|(\nabla \mathbf{u},\mathbf{u}_t,\nabla\mathbf{u}_t)\|_{L^2}^2\mm{d}s  \\
&\leq C\left(\mathcal{E}_0^2+\int_0^t\|(\varrho,\mathbf{u})\|_{L^2}^2\mm{d}s\right).
\end{aligned}\end{equation}

In addition, based on \eqref{enegetin1234}, one can infer higher regularity of $\varrho$.
Observing that $\varrho$ satisfies
\begin{equation*}
(|\varrho_{x_j}|^2)_t+\mathbf{u}\cdot\nabla(|\varrho_{x_j}|^2=-2\varrho_{x_j}(
(\nabla \bar{\rho}\cdot\mathbf{u})_{x_j}+\nabla {\varrho}\cdot
\mathbf{u}_{x_j}),\ j=1,2,3,
\end{equation*}
we integrate over $\Omega$, sum over $j$ and use (\ref{0102}) to infer that
\begin{equation*} \begin{aligned}
\frac{d}{dt}\int|\nabla
\varrho(t)|^2\mathrm{d}\mathbf{x}\leq & C\int\left(|\nabla \mathbf{u}||\nabla
\varrho|+(|\mathbf{u}|+|\nabla \mathbf{u}|\right)|\nabla\varrho|\mathrm{d}\mathbf{x}\\
\leq &C(\|\nabla \mathbf{u}\|_{W^{1,6}}\|\nabla\varrho\|_{L^2}+\|\mathbf{u}\|_{H^1})\|\nabla\varrho\|_{L^2},
\end{aligned}  \end{equation*}
whence,
 \begin{equation*} \begin{aligned}
\frac{d}{dt}\|\nabla\varrho (t)\|_{L^2}
\leq & C\big[\|\nabla \mathbf{u}\|_{W^{1,6}}\|\nabla\varrho\|_{L^2}+\|\mathbf{u}\|_{H^1}\big].\end{aligned}
\end{equation*}
Applying Gronwall's inequality, one finds that
\begin{equation}\label{L6de}\begin{aligned}
\|\nabla \varrho(t)\|^2_{L^2}\leq \left(\|\nabla\varrho_0\|_{L^2}+C\int_0^t
\|\mf{u}\|_{H^1}\mm{d}s\right)e^{C\sqrt{t}\|\nabla \mathbf{u}\|_{L^2((0,t),W^{1,6})}}.\end{aligned}
\end{equation}

On the other hand, arguing analogously to (\ref{0403}), we can deduce from (\ref{0105})$_1$ that
\begin{equation}\label{0437}
\frac{d}{dt}\|\varrho (t) \|_{L^6}\leq 6\|\bar{\rho}'\|_{L^\infty}\|\mathbf{u} (t) \|_{L^6}
\leq C \|\nabla \mathbf{u} (t) \|_{L^2}.
\end{equation}
Moreover, it follows from the classical regularity theory on the Stokes equations that
\begin{equation}\label{stokes6}
\begin{aligned}\|\nabla^2 \mathbf{u} \|_{L^6}^2+\|\nabla q\|_{L^6}^2\leq &C(\| \mathbf{u}_t\|_{L^6}^2
+ \|\mathbf{u}\cdot\nabla \mathbf{u}\|_{L^6}^2 +\|\varrho \|_{L^6}^2)  \\
\leq &C(\|\nabla \mathbf{u}_t\|_{L^2}^2
+ \|\mathbf{u}\|_{H^2}^2\| \nabla \mathbf{u}\|_{H^1}^2 +\|\varrho\|_{L^6}^2),
\end{aligned}\end{equation}
while making use of \eqref{enegetin1234}--\eqref{stokes6}, we infer that
\begin{equation}\|\nabla
\varrho(t)\|^2_{L^2}\leq
C\left(\mathcal{E}_0^2,\|\nabla\varrho_0\|_{L^2},
\sup_{0\leq s\leq t}\mathcal{E}(s),t\right),
\end{equation}
where the constant $C(\cdots\cdot)$ is nondecreasing on its four variables.

Now we mention that the local existence of strong solutions   to the 3D
nonhomogeneous incompressible Navier-Stokes equations have been
established, see \cite{CYKHOM1str,CYKHOM1ufda} for example. In
particular, by a slight modification in arguments, one can follow the proof  in \cite{CYKHOM1str} to obtain a local existence
result of a unique strong solution $(\rho,\mf{v})\in C^0([0,T^{\max}],H^1\times H^2)$ to
the perturbed problem (\ref{0105})--(\ref{0107}); moreover the strong solution satisfies the above\emph{ a
priori }estimates.  The proof is standard by means of energy
estimates, and hence we omit it here.
Consequently, summing up the above\emph{ a priori }estimates, we can
arrive at the following conclusion:
\begin{pro} \label{pro:0401} Assume that the RT density profile $\bar{\rho}$ satisfies \eqref{0102}.
For any given initial data $(\varrho_0,\mf{u}_0)\in (H^1\cap L^\infty)\times H^2$ satisfying
$\inf_{\mathbf{x}\in\Omega}\{\rho_0(\mathbf{x})\}>0$, and  the compatibility conditions
$\mf{u}|_{\partial\Omega}=\mf{0}$ and $\mm{div}\mf{u}_0=0$, there exist a $T>0$ and a unique strong solution
$(\varrho,\mf{u})\in C^0([0,T],L^2\times H^2)$ to the perturbed
problem \eqref{0105}--\eqref{0107}. Moreover, there is a constant $\bar{\delta}_0\in (0,1]$, such that
if $\mathcal{E}(t)\leq \bar{\delta}_0$ on $[0,T]$, then the strong solution satisfies
\begin{equation}\label{energyinequality}
\mathcal{E}^2(t)+\|(\mathbf{u}_{t},\nabla q)(t)\|_{L^2}^2 +\int_0^t\|(\nabla\mathbf{u},\mathbf{u}_t,\nabla
\mathbf{u}_t)(s)\|_{L^2}^2\mm{d}s \leq C\left(\mathcal{E}_0^2+\int_0^t\|(\varrho,\mathbf{u})(s)\|_{L^2}^2\mm{d}s\right),
\end{equation}
where the constant $C$ only depends on $\mu$, $g$, $\bar{\rho}$ and $\Omega$.
\end{pro}

\section{Proof of the nonlinear instability}\label{sec:05}

Now we are in a position to prove Theorem \ref{thm:0102}
by adopting the basic ideas in \cite{GYSWIC,JJTIIA,wang2011viscous}. First,
 in view of Theorem \ref{thm:0101}, we can construct a linear solution
\begin{equation}\label{0501}
\left(\varrho^\mm{l},
{\mathbf{u}}^\mm{l}\right)=e^{{\Lambda t}}
\left(\bar{\varrho}_0,
\bar{\mathbf{u}}_0\right)\in H^2\times  H^2
\end{equation}
to \eqref{0106}--\eqref{0108} with the initial data $(\bar{\varrho}_0,\bar{\mathbf{u}}_0)
\in H^2\times  H^2$ and $\mm{div}\bar{\mf{u}}_0=0$. Moreover, this solution satisfies
\begin{eqnarray}\label{n0502}
&&\|\bar{\rho}_0\|_{L^2}\|{\bar{u}}_{03}\|_{L^2}\|({\bar{u}}_{01},\bar{u}_{02})\|_{L^2}>0,\\
&&\nonumber {\mathcal{E}}((\bar{\varrho}_0,\bar{\mf{u}}_0))
=\sqrt{
\|\varrho_0\|_{L^2}^2+\|\mathbf{u}_0\|_{H^2}^2}=1,
\end{eqnarray}
where $\bar{u}_{30}$ and $u_3^\mm{l}$ stand for the third component of $\bar{\mf{u}}_0$ and $\mf{u}^{\mm{l}}$,
respectively.

Denote
$(\varrho_0^\delta,\mf{u}_0^\delta)
:=\delta (\bar{\varrho}_0,\bar{\mf{u}}_0)$, and $C_1:=\|(\bar{\varrho_0},\bar{\mf{u}}_0 )\|_{L^2}$.
Keping in mind that the condition $\inf_{\mathbf{x}\in\Omega}\{\bar{\rho}(\mathbf{x})\}>0$ and the embedding theorem
$H^2\hookrightarrow L^\infty$, we can choose a sufficiently small $\delta$, such that
$\inf_{\mathbf{x}\in\Omega}\{\rho_0^\delta(\mathbf{x}):=\varrho_0^\delta(\mathbf{x})+\bar{\rho}(\mathbf{x})\}>0$.
Hence, by virtue of Proposition \ref{pro:0401}, there is a $\tilde{\delta}\in (0,1)$,
such that for any $\delta<\tilde{\delta}$, there exists a unique local solution
$(\varrho^\delta,\mathbf{u}^\delta)\in C([0,T],H^1\times H^2)$ to (\ref{0105}), emanating
from the initial data $(\varrho_0^\delta,\mf{u}_0^\delta)$ with ${\mathcal{E}}((\varrho_0^\delta,\mf{u}_0^\delta))=\delta$.
Let $C>0$ and $\bar{\delta}_0>0$ be the same constants as in Proposition \ref{pro:0401} and
$\delta_0=\min\{\tilde{\delta},\bar{\delta}_0\}$. Let $\delta\in (0,\delta_0)$ and
 \begin{equation}\label{times}
 T^{\delta}=\frac{1}{\Lambda}\mm{ln}\frac{2\varepsilon_0}{\delta}>0,\quad\mbox{i.e.,}\;
 \delta e^{\Lambda T^\delta}=2\varepsilon_0, \end{equation}
where the value of $\varepsilon_0$, independent of $\delta$, is sufficiently
small and will be fixed later.

We then define
 \begin{equation*}
T^*=\sup\left\{t\in (0,T^{\max})\left|~{\mathcal{E}}((\varrho^\delta,
{\mathbf{u}}^\delta )(t))\leq {\delta_0}\right.\right\}\end{equation*}
 and
   \begin{equation*}
   T^{**}=\sup\left\{t\in (0,T^{\max})\left|~\left\|\left(\varrho^\delta,
{\mathbf{u}}^\delta\right)(t)\right\|_{{L}^2}\leq 2\delta C_1e^{\Lambda t}\right\}\right., \end{equation*}
where $T^{\mm{max}}$ denotes the maximal time of existence. Obviously, $T^*T^{**}>0$, and furthermore,
 \begin{eqnarray}\label{0502n1}
&&{\mathcal{E}}(\left(\varrho^\delta, {\mathbf{u}}^\delta\right)(T^*))={\delta_0}\quad\mbox{ if }T^*<\infty ,\\
\label{0502n111}  &&  \left\|\left(\varrho^\delta, {\mathbf{u}}^\delta,\right)(T^{**})\right\|_{{L}^2}
=2\delta C_1e^{\Lambda T^{**}}\quad\mbox{ if }T^{**}<T^{\max}.
\end{eqnarray}
Then for all $t\leq \min\{T^\delta,T^*,T^{**}\}$, we deduce from the estimate \eqref{energyinequality} and the
definition of $T^*$ and $T^{**}$ that
 \begin{equation}\begin{aligned}\label{0503}
{\mathcal{E}}^2\big( (\varrho^\delta, {\mathbf{u}}^\delta)(t)\big) +\|\mathbf{u}_{t}^\delta(t)\|_{L^2}^2
\leq &C \delta^2 {\mathcal{E}}^2(\left(\bar{\varrho}_0,
\bar{\mathbf{u}}_0 \right))+C\int_0^t\left\|\left(\varrho^\delta,
{\mathbf{u}}^\delta \right)(s)\right\|_{L^2}^2\mm{d}s\\
\leq & C\delta^2+4 CC_1^2\delta^2e^{2\Lambda t}/(2\Lambda)
\leq C_2\delta^2e^{2\Lambda t}
   \end{aligned}
 \end{equation}
for some constant $C_2>0$ independent of $\delta$.

Let $(\varrho^{\mathrm{d}}, {\mathbf{u}}^{\mathrm{d}})=(\varrho^{\delta},
{\mathbf{u}}^{\delta})-\delta(\varrho^{\mathrm{l}}, {\mathbf{u}}^{\mathrm{l}})$.
Noting that $(\varrho^\mm{a}_\delta,\mf{u}^{\mm{a}}_\delta):= \delta(\varrho^{\mm{l}},\mf{u}^{\mm{l}})$
is also a linear solution to \eqref{0106}--\eqref{0108} with the initial data
$(\varrho_0^\delta,\mathbf{u}_0^\delta)\in H^2\times H^2$, we find that $(\varrho^{\mathrm{d}}, {\mathbf{u}}^{\mathrm{d}})$
satisfies the following non-homogenous equation:
\begin{equation}\label{h0407}\left\{\begin{array}{ll}
  \varrho_t^{\mathrm{d}}+\bar{\rho}'{u}_3^{\mathrm{d}}= -{{\mathbf{u}}}^{\delta}\cdot \nabla\varrho^{\delta}, \\[1mm]
  \bar{\rho}\mathbf{u}_t^{\mathrm{d}}   -\mu \Delta \mathbf{u}^{\mathrm{d}}+\nabla q^{\mathrm{d}}+g\varrho^{\mathrm{d}} \mf{e}_3
=   -( \varrho^{\delta}+\bar{\rho})\mathbf{u}^{\delta}\cdot\nabla
                             \mathbf{u}^{\delta}-\varrho^\delta\mathbf{u}^{\delta}_t,\\[1mm]
 \mathrm{div}\mathbf{u}^{\mathrm{d}}={0}
 \end{array}\right.\end{equation}
 with initial data
 \begin{equation*} (\varrho^{\mathrm{d}}(0),
{\mathbf{u}}^{\mathrm{d}}(0))=
{\mathbf{0}},\quad \mathrm{div}\mf{u}_0^{\mathrm{d}}=0.
   \end{equation*}
Multiplying \eqref{h0407}$_1$ and \eqref{h0407}$_2$ by $\varrho^{\mathrm{d}}/\bar{\rho}'$ and $\mf{u}^{\mathrm{d}}/g$,
respectively, adding the resulting equalities, and integrating
over $\Omega$ (by parts), one has that for all $t\leq \min\{T^\delta,T^*,T^{**}\}$,
 \begin{equation}\label{energyford}
\begin{aligned}
& \frac{d}{dt}\int\left(\frac{|\varrho^{\mathrm{d}}|^2}{\bar{\rho}'} + \frac{\bar{\rho}|\mathbf{u}^{\mathrm{d}}|^2}{g}\right)
\mm{d}\mf{x}+\int \left(\frac{\mu|\nabla \mathbf{u}^{\mathrm{d}}|^2}{g}+2\varrho^{\mathrm{d}} u^\mm{d}_3 \right)\mm{d}\mf{x} \\
&=-\int\left(\frac{{\bar{\rho}''}{u^{\delta}_3}|\varrho^{\delta}|^2}{|\bar{\rho}'|^2}
-\frac{{\bar{\rho}''}{u^{\delta}_3}\varrho^{\delta}\varrho^{\mm{a}}_\delta}{|\bar{\rho}'|^2}
+\frac{\varrho^{\delta}{\mf{u}^{\delta}}\cdot\nabla \varrho^\mm{a}_{\delta}}{\bar{\rho}'}
\right)\mm{d}\mf{x}  \\
&\quad -\int\left(( \varrho^{\delta}+\bar{\rho})\mathbf{u}^{\delta}\cdot\nabla
 \mathbf{u}^{\delta}+\varrho^\delta\mathbf{u}^{\delta}_t\right)
\cdot\mathbf{u}^{\mathrm{d}}\mm{d}\mf{x}.
 \end{aligned}  \end{equation}

Next, we control the terms on the right hand of the above equality.
 By \eqref{0503}, \eqref{0501} and the embedding theorem $H^2\hookrightarrow L^\infty$, we obtain
   \begin{equation*}
\begin{aligned}
&\left|-\int\left(\frac{{\bar{\rho}''}{u^{\delta}_3}|\varrho^{\delta}|^2}{|\bar{\rho}'|^2}
-\frac{{\bar{\rho}''}{u^{\delta}_3}\varrho^{\delta}\varrho^{\mm{a}}_\delta}{|\bar{\rho}'|^2}
+\frac{\varrho^{\delta}{\mf{u}^{\delta}}\cdot\nabla \varrho^\mm{a}_{\delta}}{\bar{\rho}'}
\right)\mm{d}\mf{x}\right|\\
&\leq \left\|\frac{\bar{\rho}''}{\bar{\rho}'^2}
\right\|_{L^\infty}\|u^{\delta}_3\|_{L^\infty} \|\varrho^{\delta}\|_{L^2}\left(\|\varrho^{\delta}\|_{L^2}
+\|\varrho^{\mm{a}}_{\delta}\|_{L^2}\right)
+\left\|\frac{1}{\bar{\rho}'} \right\|_{L^\infty}\|\mf{u}^{\delta}\|_{L^\infty}
\|\varrho^{\delta}\|_{L^2}\|\nabla\varrho^{\mm{a}}_{\delta}\|_{L^2}
\\
&\leq C_0\left(\left\|\frac{\bar{\rho}''}{\bar{\rho}'^2}
\right\|_{L^\infty}  +\left\|\frac{1}{\bar{\rho}'}
\right\|_{L^\infty}\right) {\mathcal{E}}^2(\left(\varrho^\delta,
{\mathbf{u}}^\delta\right))\left[{\mathcal{E}}(\left(\varrho^\delta, {\mathbf{u}}^\delta\right))
+\delta e^{\Lambda t}(\|\bar{\varrho}_0\|_{L^2} +\|\nabla\bar{\varrho}_0\|_{L^2})\right]
\\
& \leq C_3\delta^3e^{3\Lambda t}
 \end{aligned}  \end{equation*}
 for some constant $C_3$ dependent of $\bar{\rho}$, $\|\bar{\varrho}_0\|_{H^1}$,
 $C_2$ and the embedding constant $C_0$. Similarly, one has
  \begin{equation*}
\begin{aligned}
&\left|-\int\left(( \varrho^{\delta}+\bar{\rho})\mathbf{u}^{\delta}\cdot\nabla
                             \mathbf{u}^{\delta}+\varrho^{\delta}\mathbf{u}^{\delta}_t\right)
                             \cdot\mathbf{u}^{\mathrm{d}}\mm{d}\mf{x}\right|
                            \\
& \leq  \left(\| \varrho^{\delta}+\bar{\rho}\|_{L^\infty}
\|\mathbf{u}^{\delta}\|_{L^\infty}\|\nabla\mathbf{u}^{\delta}\|_{L^2}
\|(\mathbf{u}^{{\delta}}-\mathbf{u}^{\mm{a}}_\delta)\|_{L^2}
+\|\varrho^\delta\|_{L^2}\|\mathbf{u}^{\delta}_t \|_{L^2}\|(\mathbf{u}^{{\delta}}-\mathbf{u}^{\mm{a}}_\delta)\|_{L^\infty}\right)
\\
& \leq  C_4\delta^3e^{3\Lambda t},
 \end{aligned}  \end{equation*}
 where the constant $C_4$ depends on $\bar{\rho}$, $\|(\bar{\varrho}_0,\bar{\mf{u}}_0)\|_{H^2}$, $C_0$ and $C_2$.

On the other hand, similarly to \eqref{supress}, one deduces that
 \begin{equation*}
\begin{aligned}-\int \left(\frac{\mu|\nabla \mathbf{u}^{\mathrm{d}}|^2}{g}+2\varrho^{\mathrm{d}} u^\mm{d}_3 \right)\mm{d}\mf{x}\leq
\Lambda\int\left(\frac{|\varrho^{\mathrm{d}}|^2}{\bar{\rho}'} + \frac{\bar{\rho}|\mathbf{u}^{\mathrm{d}}|^2}{g}\right)\mm{d}\mf{x}
 \end{aligned}  \end{equation*}
 Thus, substituting the above three estimates into \eqref{energyford}, we conclude that
   \begin{equation*}
\begin{aligned}
& \frac{d}{dt}\int\left(\frac{|\varrho^{\mathrm{d}}|^2}{\bar{\rho}'} + \frac{\bar{\rho}|\mathbf{u}^{\mathrm{d}}|^2}{g}\right)\mm{d}\mf{x}
+\int \left(\frac{\mu|\nabla \mathbf{u}^{\mathrm{d}}|^2}{g}+2\varrho^{\mathrm{d}} u^\mm{d}_3 \right)\mm{d}\mf{x}   \\
&\leq \Lambda\int\left(\frac{|\varrho^{\mathrm{d}}|^2}{\bar{\rho}'} + \frac{\bar{\rho}|\mathbf{u}^{\mathrm{d}}|^2}{g}\right)\mm{d}\mf{x}
+(C_3+C_4)\delta^3e^{3\Lambda t}    \end{aligned}  \end{equation*}
 Applying Gronwall's inequality to the above inequality, one obtains that for all $t\leq \min\{T^\delta,T^*,T^{**}\}$,
 \begin{equation*}
\begin{aligned}
 \int\left(\frac{|\varrho^{\mathrm{d}}|^2}{\bar{\rho}'} + \frac{\bar{\rho}|\mathbf{u}^{\mathrm{d}}|^2}{g}\right)\mm{d}\mf{x}
&\leq \int_0^t (C_3+C_3)\delta^3e^{3\Lambda s}e^{\Lambda(t-s)}\mm{d}s\\
&\leq (C_3+C_4)\delta^3e^{3\Lambda t}/2\Lambda,
 \end{aligned}  \end{equation*}
 which yields
 \begin{equation}\label{efd}  \begin{aligned}
 \|(\varrho^{\mathrm{d}},\mathbf{u}^{\mathrm{d}}  )\|_{L^2}^2
\leq C_5\delta^3e^{3\Lambda t}\quad\mbox{ for some constant } C_5.
 \end{aligned}  \end{equation}

 Now, we claim that
\begin{equation}\label{n0508}
T^\delta=\min\left\{T^\delta,T^*,T^{**}\right\},
 \end{equation}
provided that small $\varepsilon_0$ is taken to be
 \begin{equation}\label{defined}
\varepsilon_0=\min\left\{\frac{{\delta_0}}
{4\sqrt{C_2}},\frac{C_1^2}{8C_5},\frac{m_0^2}{C_5} \right\},
 \end{equation}where
$m_0= \min\{
 \|\bar{\varrho}_0\|_{L^2},\|\bar{u}_{03}\|_{L^2},
\|(\bar{u}_{01},\bar{u}_{02})\|_{L^2}\}>0$ due to \eqref{n0502}.

 Indeed, if $T^*=\min\{T^{\delta},T^*,T^{**}\}$, then
$T^*<\infty$. Moreover, from \eqref{0503} and \eqref{times} we get
 \begin{equation*}
{\mathcal{E}}(\left(\varrho^\delta,
{\mathbf{u}}^\delta
\right)(T^*))\leq \sqrt{C_2}\delta e^{\Lambda T^*}
\leq \sqrt{C_2}\delta e^{\Lambda T^\delta}=2\sqrt{C_2}\varepsilon_0<{\delta_0},
 \end{equation*}
 which contradicts with \eqref{0502n1}. On the other hand, if $T^{**}=\min\{T^{\delta},T^*,T^{**}\}$, then $T^{**}<T^{\mm{max}}$.
 Moreover, in view of \eqref{0501}, \eqref{times} and \eqref{efd}, we see that
 \begin{equation*}\begin{aligned}
 \left\|\left(\varrho^\delta,
{\mathbf{u}}^\delta
\right)(T^{**})\right\|_{L^2}
\leq  & \left\|\left(\varrho^\mm{a}_{\delta},
{\mathbf{u}}^\mm{a}_{\delta}
\right)(T^{**})\right\|_{L^2} +\left\|\left(\varrho^{\mathrm{d}},
{\mathbf{u}}^{\mathrm{d}}
\right)(T^{**})\right\|_{L^2} \\
\leq  &\delta \left\|\left(\varrho^\mm{l},
{\mathbf{u}}^{\mm{l}}
\right)(T^{**})\right\|_{L^2}+\sqrt{C_5}\delta^{3/2}e^{3\Lambda T^{**}/2} \\
\leq & \delta C_1e^{\Lambda T^{**}}+\sqrt{C_5}\delta^{3/2} e^{3\Lambda T^{**}/2}
\leq \delta e^{\Lambda T^{**}}(C_1+\sqrt{2C_5\varepsilon_0})\\
<&2\delta C_1  e^{\Lambda T^{**}},
 \end{aligned} \end{equation*}
which also contradicts with \eqref{0502n111}. Therefore, \eqref{n0508} holds.

 Finally, we again use \eqref{defined} and \eqref{efd} to deduce that
 \begin{equation*}\begin{aligned}
 \|\varrho^{\delta}(T^\delta)\|_{L^2}\geq &
\|\varrho^{\mathrm{a}}_{\delta}(T^{\delta})\|_{L^2}-\|\varrho^{\mm{d}}(T^{\delta})\|_{L^2}
= \delta\|\varrho^{\mathrm{l}}(T^{\delta})\|_{L^2}-\|\varrho^{\mm{d}}(T^{\delta})\|_{L^2} \\
 \geq & \delta e^{\Lambda T^\delta}\|\bar{\varrho}_{0}\|_{L^2} -\sqrt{C_5}\delta^{3/2}e^{3\Lambda^* T^{\delta}/2}
 \\
 \geq & 2\varepsilon_0\|\bar{\varrho}_{0}\|_{L^2} -\sqrt{C_5}\varepsilon_0^{3/2}
 \geq 2m_0\varepsilon_0 -\sqrt{C_5}\varepsilon_0^{3/2} \geq m_0\varepsilon_0,
 \end{aligned}      \end{equation*}
 Similar, we also have \begin{equation*}\begin{aligned}
 \|u_3^{\delta}(T^\delta)\|_{L^2}
 \geq  2m_0\varepsilon_0 -\sqrt{C_5}\varepsilon_0^{3/2}\geq m_0\varepsilon_0,
 \end{aligned}      \end{equation*}
and
 \begin{equation*}\begin{aligned}
 \|(u_1^{\delta},u_2^{\delta})(T^\delta)\|_{L^2}\geq
   2m_0\varepsilon_0 -\sqrt{C_5}\varepsilon_0^{3/2}
\geq m_0\varepsilon_0,
 \end{aligned}      \end{equation*}
 where $u^{\delta}_{i}(T^{\delta})$ denote the $i$-th component of
$\mf{u}^{\delta}(T^{\delta})$ for $i=1$, $2$, $3$.
This completes the proof of Theorem \ref{thm:0102} by defining $\varepsilon :=m_0\varepsilon_0$.


\section{Proof of the stability}\label{sec:06}
Similar to the proof of Proposition \ref{pro:0401}, we see that to get Theorem \ref{thm:0103}, it suffices to deduce
the \emph{a priori } estimates \eqref{nnn0117}--\eqref{nnn0118}. In what follows, we denote by $C$ a
generic positive constant which may depend on $\mu$, $g$, $\bar{\rho}$ and $\Omega$.

Let $(\varrho,\mf{u})$ solve the linearized problem \eqref{0106}--\eqref{0108} with an associated pressure $q$.
Recalling the condition $\sup_{\mf{x}\in \Omega}\bar{\rho}'(\mf{x})<0$ in Theorem \ref{thm:0103},
we can derive from the linearized mass and momentum equations that
\begin{equation}\label{mass6}\frac{d}{dt}
\int\frac{g{\varrho^2(t)}}{{-\bar{\rho}'}}\mm{d}\mf{x}=2 g\int\varrho
u_3\mathrm{d}\mathbf{x}
\end{equation}
and
\begin{equation}\label{energyequality}\frac{d}{dt}\int
\bar{\rho}|\mathbf{u}|^2(t)\mathrm{d}\mathbf{x}+2\mu\int|\nabla \mathbf{u}|^2
\mathrm{d}\mathbf{x}= -2g\int {\varrho }{u}_3\mathrm{d}\mathbf{x}.\end{equation}
Adding the above two equalities, one gets
\begin{equation*}
\begin{aligned}&\frac{d}{dt}
\left\|\left(\sqrt{\frac{g}{-\bar{\rho}'}}\varrho,
\sqrt{\bar{\rho}}\mathbf{u}\right)(t)\right\|^2_{L^2}+2\mu\|\nabla \mathbf{u}\|^2_{L^2}=0,
\end{aligned}\end{equation*}
which immediately yields the estimate \eqref{nnn0117}, i.e,
  \begin{equation}\label{n0601}
\|(\varrho, \mathbf{u})(t)\|^2_{L^2}+\int_0^t\|\nabla
 \mathbf{u}\|^2_{L^2}\mm{d}s\leq C\|(\varrho_0,\mathbf{u}_0)\|^2_{L^2}.
\end{equation}

Analogous to \eqref{nnn0314}, we find that $\mf{u}$ satisfies
\begin{equation*}
\frac{d}{dt}\int \left(\bar{\rho}|\mathbf{u}_t|^2 -g\bar{\rho}'{u}_3^2\right)\mathrm{d}\mathbf{x}+
2\mu\int |\nabla\mathbf{u}_t|^2\mm{d}\mf{x}=0,
\end{equation*}
whence
\begin{equation*}
\|({u}_3,\mathbf{u}_t)\|^2_{L^2}+ \int_0^t \|\nabla\mathbf{u}_t(s)\|^2_{L^2}\mm{d}s
\leq  C\|(\varrho_0,{u}_3,\Delta\mf{u}_0)\|_{L^2}^2.
\end{equation*}
On the other hand, by the classical regularity theory on the Stokes equations, one infers that
\begin{equation*}
\|(\nabla^2 \mathbf{u},\nabla q)\|_{L^2}^2\leq
C\|-\bar{\rho} \mathbf{u}_t-\varrho {g}e_3\|_{L^2}^2\leq  C \| (\varrho,\mathbf{u}_t)\|_{L^2}^2,
\end{equation*}
which combined with \eqref{n0601} implies
  \begin{equation}\label{123456}
\|( \nabla^2\mf{u},\mf{u}_t,\nabla q)(t)\|^2_{L^2}+\int_0^t\|\nabla\mathbf{u}_t(s)\|^2_{L^2}\mm{d}s
\leq\|(\varrho_0, \mathbf{u}_0,\Delta\mf{u}_0)\|^2_{L^2}.
\end{equation}

Finally, by the well-known Gagliardo-Nirenberg interpolation inequality (see \cite[Theorem]{NLOE}), we have
 \begin{eqnarray}\label{nireberg}
  \|\nabla \mf{u}\|_{L^{2}}\leq C( \|\mf{u}\|_{L^2}^{\frac{1}{2}}\|\nabla^{2}\mf{u}\|_{L^{2}}^{\frac{1}{2}}
 +\| \mf{u}\|_{L^{2}}), \end{eqnarray}
 which yields that
 \begin{eqnarray}\label{gradiane}
  \|\nabla\mf{u}\|_{L^{2}}\leq C \|(\mf{u},\nabla^{2}\mf{u})\|_{L^{2}}.
\end{eqnarray}
Consequently, we deduce from \eqref{n0601} and \eqref{123456} that
\begin{equation}\label{n0604}
\|(\mf{u}_t,\nabla q)(t)\|^2_{L^2}+\|\nabla\mf{u}(t)\|^2_{H^1}+\int_0^t\|\nabla
\mathbf{u}_s\|^2_{L^2}\mm{d}s\leq\|(\varrho_0,
\mathbf{u}_0,\Delta\mf{u}_0)\|^2_{L^2},\end{equation}
which yields \eqref{nnn0117}.

By \cite[Theorem 1.68]{NASII04}, Sobolev's inequality, \eqref{n0601} and \eqref{n0604}, we have
\begin{equation}\label{convergence}
\begin{aligned}\int_0^t\left|\frac{d}{dt}\|(\mf{u},\nabla\mf{u})(t)\|_{L^2}^2\right|\mm{d}s
 &\leq 2\int_0^t\left(\|\mf{u}\|_{L^2}\|\mf{u}_t\|_{H^{-1}} + \|\nabla\mf{u}\|_{L^2}\|\nabla\mf{u}_t\|_{H^{-1}}\right) \\
&\leq C\int_0^t \|(\mf{u},\nabla\mf{u},\mf{u}_t,\nabla\mf{u}_t)\|_{L^2}^2\mm{d}s  \\
&\leq C\int_0^t \|(\nabla \mf{u},\nabla \mf{u}_t)\|_{L^2}^2\mm{d}s\leq C\|(\varrho_0,\mathbf{u}_0,\Delta\mf{u}_0)\|^2_{L^2}.
\end{aligned}\end{equation}
Hence $\|\mf{u}(t)\|_{H^1}\in W^{1,1}(0,\infty)$, which implies
\eqref{nnn0119}, i.e.,  $\|\mathbf{u}(t)\|_{H^1}\to 0$ as $t\to\infty$.

In addition, for any $T>0$, using \eqref{n0601} and \eqref{n0604}, we can derive
more \emph{a priori} estimates $\nabla\rho\in L^\infty((0,T),L^2(\Omega))$,
$\rho_t\in L^\infty((0,T),H^1(\Omega))$ and $\nabla^3\mf{u}\in L^2((0,T),L^2(\Omega))$
for the linearized equations \eqref{0108}. With these \emph{a priori} estimates and \eqref{nnn0119},
adapting the Faedo-Galerkin approximation scheme as done in \cite{CYKHOM1str}, we can easily obtain the
first assertion in Theorem \ref{thm:0103}.

To show the \emph{a priori} estimates of the second assertion in Theorem \ref{thm:0103}, we let $(\varrho,\mf{u})$
solve the perturbed problem \eqref{0105}--\eqref{0107} with the associated pressure $q$.
Recalling that $\bar{\rho}'$ is constant, we deduce from \eqref{0105} that $(\varrho,\mf{u})$ satisfies
\begin{equation}\label{weilon}
\begin{aligned}&\frac{d}{dt}
\left\|\left(\sqrt{\frac{g}{-\bar{\rho}'}}\varrho,
\sqrt{{\rho}}\mathbf{u}\right)(t)\right\|^2_{L^2}+2\mu\|\nabla
 \mathbf{u}\|^2_{L^2}=0,
\end{aligned}\end{equation}
where $\rho=\varrho +\bar{\rho}$. Therefore, we have
  \begin{equation}\label{n060111}
\|(\varrho,\mathbf{u})(t)\|^2_{L^2}+\int_0^t\|\nabla\mathbf{u}\|^2_{L^2}\mm{d}s\leq C\|(\varrho_0,\mathbf{u}_0)\|^2_{L^2}.
\end{equation}

If we argue in a manner similar to that in the derivation of \eqref{0425}, we find that
\begin{equation*} \begin{aligned}
&\frac{ 1}{2}\frac{d}{dt}\int\left(\rho|\mathbf{u}_{t}|^2-g\bar{\rho}'u_3^2\right)\mm{d}\mf{x} +\mu|\nabla\mathbf{u}_t|^2\mm{d}\mf{x}\\
&\leq \int2\rho|\mathbf{u}||\mathbf{u}_t||\nabla\mathbf{u}_t|+\rho|\mathbf{u}||\mathbf{u}_t||\nabla\mathbf{u}|^2
+\rho |\mathbf{u}|^2|\mathbf{u}_t||\nabla^2 \mathbf{u}|  \\
&\quad +\rho|\mathbf{u} |^2|\nabla \mathbf{u}||\nabla\mathbf{u}_t|+\rho|\mathbf{u}_t|^2|\nabla
\mathbf{u}|+g\varrho|\mathbf{u}||\nabla\mathbf{u}_t|:=\sum_{i=1}^6K_j.\end{aligned}
\end{equation*}
Letting
$\sup_{\mathbf{x}\in\Omega}\{\rho_0({\mathbf{x}})\}\leq K$, since
\begin{equation*}
0<\inf_{\mathbf{x}\in\Omega}\{\rho_0(\mathbf{x})\} \leq \rho(t,\mf{x})\leq \sup_{\mathbf{x}
\in\Omega}\{\rho_0({\mathbf{x}})\}\quad\mbox{ for any }t>0\mbox{ and }\mf{x}\in \Omega ,
\end{equation*}
we have
\begin{equation*}
\|\varrho(t)\|_{L^\infty}\leq K+\|\bar{\rho}\|_{L^\infty}.
\end{equation*}
Hence,
\begin{eqnarray*}
K_6\leq C(\varepsilon)\|\varrho\|_{L^3}^2
\|\mathbf{u}\|_{L^6}^2+\varepsilon\|\nabla\mathbf{u}_t\|_{L^2}^2
\leq C(\varepsilon)(K+\|\bar{\rho}\|_{L^\infty})^{\frac{2}{3}}\|\varrho\|_{L^2}^{\frac{4}{3}}
\|\nabla \mf{u}\|_{L^2}^2+\varepsilon\|\nabla\mathbf{u}_t\|_{L^2}^2.
\end{eqnarray*}

Similar to \eqref{jfjs0427}, we conclude
\begin{equation}\label{sxk0427}
\begin{aligned}
&\frac{d}{dt}\left\|\left(\sqrt{\rho}\mathbf{u}_{t},\sqrt{-g\bar{\rho}'}u_3\right)(t)
\right\|_{L^2}^2+\mu\|\nabla \mathbf{u}_t\|_{L^2}^2 \\
&\leq  C(K^{\frac{2}{3}}+1)\Big[\|\nabla \mathbf{u}\|^4_{L^2}(\|\mathbf{u}_t\|^2_{L^2}
+\|\nabla \mathbf{u}\|_{H^1}^2)+\|\varrho\|_{L^2}^\frac{4}{3}\|  \nabla \mf{u}\|_{L^2}^2\Big],
\end{aligned}\end{equation}
while in a manner similar to \eqref{0433}, one obtains
\begin{equation*}
\begin{aligned}\|(\nabla^2 \mathbf{u} ,\nabla q)\|_{L^2}^2
\leq C[\|(\varrho,\mathbf{u}_t)\|_{L^2}^2+\|\nabla \mathbf{u}\|_{L^2}^6]+\frac{1}{2}\|\nabla \mathbf{u}\|_{H^1}^2,
\end{aligned}\end{equation*}
which, together with \eqref{nireberg}, yields
\begin{equation}\label{984}
\begin{aligned}\|(\nabla^2 \mathbf{u} ,\nabla q)\|_{L^2}^2
\leq C\|(\varrho,\mathbf{u},\mathbf{u}_t)\|_{L^2}^2,
\end{aligned}\end{equation}
provided that $\|\nabla \mathbf{u}\|_{L^2}$ is sufficiently small.
Thus, we can make use of \eqref{weilon}, \eqref{sxk0427}, \eqref{984} and Poincar\'e's inequality to conclude that
\begin{equation*}
\begin{aligned}&\frac{d}{dt}\left\|\left(\sqrt{\frac{g}{{-\bar{\rho}'}}}\varrho,
\sqrt{\rho}\mathbf{u},\sqrt{-g\bar{\rho}'}u_3,\sqrt{\rho}\mathbf{u}_{t}\right)(t)\right\|_{L^2}^2
+\mu\|\nabla (\sqrt{2}\mathbf{u},\mathbf{u}_t)\|^2_{L^2} \\
&\leq C(K^{\frac{2}{3}}+1)\Big[\|\nabla\mathbf{u}\|^4_{L^2}\|(\varrho,\nabla\mathbf{u}_t,\nabla
\mathbf{u})\|_{L^2}^2+\|\varrho\|_{L^2}^\frac{4}{3}\|\nabla \mf{u}\|_{L^2}^2\Big].
\end{aligned}\end{equation*}
Hence, if $\|(\varrho,\nabla \mathbf{u})\|_{L^2}$, dependent of $K$, is sufficiently small, then
\begin{equation*}
\begin{aligned}\frac{d}{dt}
\left\|\left(\sqrt{\frac{g}{{-\bar{\rho}'}}}\varrho,
\sqrt{\rho}\mathbf{u},\sqrt{-g\bar{\rho}'}u_3,\sqrt{\rho}\mathbf{u}_{t}\right)(t)\right\|_{L^2}^2
+\frac{\mu}{2}\|\nabla (\sqrt{2}\mathbf{u},\mathbf{u}_t)\|^2_{L^2}\leq 0,
\end{aligned}\end{equation*}
from which it follows that
\begin{equation}\begin{aligned}\label{nonegertn}
\|(\varrho,\mathbf{u},\mathbf{u}_{t})(t)\|_{L^2}^2 +\int_0^t\|\nabla (\mathbf{u},\mathbf{u}_t)\|^2_{L^2}\mm{d}s
\leq C\|(\varrho_0, \mathbf{u}_0,\Delta\mf{u}_0)\|^2_{L^2}.
\end{aligned}\end{equation}
Consequently, using \eqref{gradiane}, \eqref{984} and \eqref{nonegertn}, we arrive at
\begin{equation}\label{061011}
\|(\mf{u}_t,\nabla q)\|^2_{L^2}+\|\nabla\mathbf{u}\|_{H^1}+\int_0^t\|\nabla(
 \mathbf{u},\mathbf{u}_t)\|^2_{L^2}\leq C\|(\varrho_0,\mathbf{u}_0,\Delta\mf{u}_0)\|^2_{L^2},
\end{equation}
provided $\|(\varrho,\nabla \mathbf{u})\|_{H^1}$, dependent of $K$, is sufficiently small.
So we shall impose smallness condition (dependent of $K$) on the initial data
 $\|(\varrho_0,\mathbf{u}_0,\Delta\mf{u}_0)\|^2_{L^2}$ to get \eqref{061011}.

In addition, we also have the estimate \eqref{convergence}, thus $\|\mathbf{u}(t)\|_{H^1}\to 0$ as $t\to\infty$. Moreover,
 for any $T>0$, we can use \eqref{n060111} and \eqref{061011} to get
more \emph{a priori} estimates $\nabla\rho\in L^\infty((0,T),L^2(\Omega))$,
$\rho_t\in L^\infty((0,T),H^1(\Omega))$ and $\nabla^3\mf{u}\in L^2((0,T),L^2(\Omega))$
for the the perturbed problem \eqref{0105}--\eqref{0107}. With the help of these \emph{a priori} estimates, we can adapt
the Faedo-Galerkin approximation as in \cite{CYKHOM1str} to easily obtain the
second assertion in Theorem \ref{thm:0103}. This completes the proof of Theorem \ref{thm:0103}.

\vspace{4mm} \noindent\textbf{Acknowledgements.}
The research of Fei Jiang was supported by NSFC (Grant Nos. 11101044, 11271051 and
11301083), and the research of Song Jiang by the National Basic Research Program
under the Grant 2011CB309705 and NSFC (Grant Nos. 11229101, 11371065).
\renewcommand\refname{References}
\renewenvironment{thebibliography}[1]{%
\section*{\refname}
\list{{\arabic{enumi}}}{\def\makelabel##1{\hss{##1}}\topsep=0mm
\parsep=0mm
\partopsep=0mm\itemsep=0mm
\labelsep=1ex\itemindent=0mm
\settowidth\labelwidth{\small[#1]}%
\leftmargin\labelwidth \advance\leftmargin\labelsep
\advance\leftmargin -\itemindent
\usecounter{enumi}}\small
\def\newblock{\ }
\sloppy\clubpenalty4000\widowpenalty4000
\sfcode`\.=1000\relax}{\endlist}
\bibliographystyle{model1b-num-names}

\end{document}